\documentclass[11pt]{amsart}
\usepackage{enumerate,amssymb,amsmath,mathrsfs}
\usepackage{amssymb,amsfonts,amsthm,amscd}
\usepackage[mathscr]{eucal}     
\usepackage{graphicx}
\usepackage[arrow, matrix, curve]{xy}
\usepackage{color}
\usepackage[margin=3cm]{geometry}
\definecolor{gray}{gray}{.75}
\definecolor{gray2}{gray}{.50}
\usepackage{tikz}

\newcommand{\dom}{\mathcal{D}}

\newcommand{\N}{\mathbb{N}}
\newcommand{\R}{\mathbb{R}}

\newcommand{\FD}{\mathrm{FD}}
\newcommand{\AT}{\mathrm{AT}}
\newcommand{\rh}{\mathrm{rh}}

\newcommand{\ka}{\kappa}
\newcommand{\ha}{\mathfrak{h}}

\newcommand{\eh}{\mathrm{e-h}}
\newcommand{\tr}{\mathrm{tr}\,}
\newcommand{\Tr}{\mathrm{Tr}\,}

\newcommand{\A}{\mathrm{\alpha}}
\newcommand{\be}{\mathrm{\beta}}
\newcommand{\w}{\mathrm{\omega}}

\newcommand{\phg}{\mathrm{phg}}

\newcommand{\e}{\epsilon}
\newcommand{\del}{\partial}
\newcommand{\calA}{\mathcal A}
\newcommand{\calU}{\mathcal U}
\newcommand{\calE}{\mathcal E}

\newcommand{\calC}{\mathcal C}
\newcommand{\calO}{\mathcal O}
\newcommand{\calV}{\mathcal V}
\newcommand{\calB}{\mathcal B}
\newcommand{\calR}{\mathcal R}
\newcommand{\calL}{\mathcal L}
\newcommand{\calD}{\mathcal D}

\newcommand{\calG}{\mathcal G}
\newcommand{\RR}{\mathbb R}
\newcommand{\wt}{\widetilde}

\newcommand{\wx}{\widetilde{x}}
\newcommand{\wy}{\widetilde{y}}
\newcommand{\wz}{\widetilde{z}}

\newcommand{\ff}{\mathrm{ff}}
\newcommand{\lf}{\mathrm{lf}}
\newcommand{\rf}{\mathrm{rf}}
\newcommand{\tf}{\mathrm{tf}}
\newcommand{\td}{\mathrm{td}}

\newtheorem{thm}{Theorem}[section]

\newtheorem{remark}[thm]{Remark}
\newtheorem{prop}[thm]{Proposition}
\newtheorem{lemma}[thm]{Lemma}
\newtheorem{defn}[thm]{Definition}

\numberwithin{equation}{section}

\begin{document}

\title{Analytic torsion on manifolds with edges}

\author{Rafe Mazzeo}
\address{Stanford University \\ Department of Mathematics \\ Stanford, CA 94305-2125}
\email{mazzeo@math.stanford.edu}

\author{Boris Vertman}
\address{University Bonn \\ Department of Mathematics \\
Endenicher Allee 60, 53115 Bonn}
\email{vertman@math.uni-bonn.de}


\date{This document was compiled on: \today}

\begin{abstract}
Let $(M,g)$ be an odd-dimensional incomplete compact Riemannian singular space with a simple edge singularity.  We study the 
analytic torsion on $M$, and in particular consider how it depends on the metric $g$. If $g$ is an admissible edge metric, we prove 
that the torsion zeta function is holomorphic near $s=0$, hence the torsion is well-defined, but possibly depends on $g$. In general 
dimensions, we prove that the torsion depends only on the asymptotic structure of $g$ near the singular stratum of $M$; when 
the dimension of the edge is odd, we prove that the analytic torsion is independent of the choice of even admissible edge metrics.
The main tool is the construction, via the methodology of geometric microlocal analysis, of the heat kernel for the Friedrichs
extension of the Hodge Laplacian in all degrees. In this way we obtain detailed asymptotics of this heat kernel and its trace.
\end{abstract}

\maketitle

\pagestyle{myheadings}
\markboth{\textsc{Analytic torsion on manifolds with edges}}{\textsc{Rafe Mazzeo and Boris Vertman}} 

\section{Introduction}
One of the key achievements in modern spectral geometry is the proof by Cheeger and M\"uller of the Ray-Singer conjecture,
which equates the analytic and Reidemeister torsions of a compact smooth odd-dimensional manifold (equipped with a flat 
Hermitian vector bundle). Since one of these quantities is analytic and the other combinatorial, their equality has many important 
applications in fields ranging from topology and number theory to mathematical physics. The original definition of analytic torsion,
and its conjectured relationship with Reidemeister torsion, appeared in the famous 1971 paper of Ray and Singer \cite{RS}.
The original proofs by Cheeger \cite{Ch-AT} and M\"uller \cite{Mue-AT} are still of great interest, but there are now 
several other proofs of this result as well, each with its own significance and leading to further generalizations. Amongst
these we mention in particular the ones based on Witten deformation \cite{Zh} and `analytic surgery' \cite{Ha}. 

Singular spaces arise naturally in many parts of mathematics, and the development of analytic techniques to study partial differential 
equations on them is a central challenge in modern geometric analysis.  Important examples of singular spaces include algebraic 
varieties and various moduli spaces, and they also appear naturally as compactifications of smooth spaces 
or as limits of families of smooth spaces under controlled degeneration.  

A natural and still-outstanding open problem is to determine whether the Cheeger-M\"uller theorem has any analogue
for compact stratified pseudomanifolds (which are perhaps the best behaved type of singular space).
We refer to \S 2 of \cite{ALMP} for a detailed explanation of the differential topological structure of these spaces.
There are many difficulties even to formulate a precise conjecture. 
The natural homology theories in this setting are the intersection homology spaces of Goresky and MacPherson; 
however, there are many of these, none necessarily preferred over the others, which leads to an ambiguity in what
one should mean by intersection Reidemeister torsion. We refer to recent work by Dai and Huang \cite{DH} for a study
of this issue for spaces with isolated conic singularities, which is already not so straightforward.  

It is not so easy to define analytic torsion in this setting either. The development of techniques to study elliptic differential
operators, in particular geometrically natural ones such as Dirac- or Laplace-type operators, on singular spaces goes back 
many years, starting from the work of Kondratiev and his school in the early 1960's. Cheeger \cite{Ch1}, \cite{Ch2} was the
first to understand the tractability of studying spectral geometry for smoothly stratified spaces endowed with `incomplete 
iterated edge metrics'.  This was the beginning of many further developments by various authors.  The analysis of partial
differential operators for spaces with conic (and the closely related asymptotically cylindrical) metrics was also studied
by Lesch \cite{Lesch}, Melrose \cite{Mel-APS}, Br\"uning and Seeley \cite{BS}, Gil and Mendoza \cite{GM} and Mooers \cite{Moo}, 
to name just a few.  Extensions to spaces with simple edge singularities were developed by the first author 
\cite{Maz-edge}, and Schulze and his collaborators \cite{Sch1}; see also Br\"uning-Seeley \cite{BS3} and more recently 
\cite{GKM}. Further extensions to spaces with 
iterated edge singularities are still in a less refined state of development, though see \cite{ALMP}.

The present paper focuses on the definition of analytic torsion on the next simplest class of spaces beyond the ones with isolated conic 
singularities. Namely, we consider spaces $(M,g)$ with a simple edge singularity and incomplete edge metric. This means the 
following. We assume that $M$ is a compact stratified space, with a single top-dimensional stratum and only one other lower
dimensional stratum $B$, which is therefore a smooth closed manifold. It is convenient, and not too misleading, to refer to the main 
stratum as $M$, and to write $\overline{M}$ when we wish to emphasize that we are talking about the entire space. We set $b = \dim B$ 
and $m = \dim M$. The stratification hypothesis provides a neighbourhood $\calU$ of $B$ which is the total space of a smooth 
bundle over $B$ with fibre $C(F)$, an open truncated cone over a compact smooth manifold $F$ of dimension $f$, i.e. \ 
$C(F) = [0,1) \times F/\sim$ with $(0,z) \sim (0,z')$ for any $z, z' \in F$. The metric $g$ on 
$M$ is arbitrary away from $\calU$, but in this neighbourhood takes the special form
\begin{equation}
\label{eq:edgemetric}
g  := dx^2 + x^2 \ka(x) + \phi^*\ha(x),
\end{equation}
where $x$ is a smooth function on $\calU \setminus B$ which restricts to a radial function on each conic fibre, $\phi: \{x = \mbox{const.}\} 
\to B$ is the fibration of each level set, and $\ka(x)$ is a family
of smooth metrics on $F$ depending smoothly on the parameter $x$ for $0 \leq x < 1$, and $\ha(x)$ is a family of smooth metrics 
on $B$ which likewise depends smoothly on $x \in [0,1)$. We give a more invariant definition of this class of metrics below. 
For simplicity we call $(M,g)$ a simple edge space.

Our goal is to study the analytic torsion of simple edge spaces. Ideally, one hopes to prove that the analytic torsion is well-defined 
and independent of the metric $g$, and that it defines an invariant which can be computed in terms of combinatorial (and perhaps
other) data. We accomplish the first part of this here: namely we prove the existence of the analytic torsion for simple edge spaces 
and show its invariance properties under various conditions.

To state our main result, let us recall some definitions and standard terminology.  Let $\Delta_{k,g}$ be the Hodge Laplace operator 
on a compact smooth manifold $M$ with respect to the metric $g$, acting on $k$-forms, and $\det'$ its zeta-regularized determinant, 
with the zero modes removed. Define the (scalar) analytic torsion $T(M,g)$ by 
\[
\log T(M,g) = \sum_{k=0}^m (-1)^k k \log \left( \det\nolimits' \Delta_{k,g} \right).
\]
The determinant line on $M$ is defined in terms of the de Rham cohomology by
\begin{align}
\det \mathcal{H}^*(M) :=\bigotimes _{k=0}^{m}\left(\bigwedge\nolimits^{\textup{top}} H^k_{dR}(M)\right)^{(-1)^{k}}.
\end{align}
Identifying $H^k_{\mathrm{dR}}(M) \cong \ker \Delta_{k,g}$, then the determinant line on $M$ inherits a natural $L^2$ Hermitian 
structure $\| \cdot \|_{L^2}$ induced from the inclusion $\ker \Delta_{k,g} \subset L^2\Omega^k(M, dV_g)$. There is
another Hermitian structure on this same line bundle, described by the norm
\[
\big\|\cdot \big\|^{RS}_{(M,g)}:=T(M, g)\big\|\cdot \big\|_{L^2};
\]
this is called either the analytic torsion, the Ray-Singer or the Quillen metric. 

Now return to the setting where $M$ has a simple edge. We introduce in \S 2 various classes of metrics on $M$. The most 
restrictive is the class of rigid (incomplete edge) metrics, which are exactly warped product conic on each conic fibre of the 
tubular neighbourhood $\calU$, and amongst these we also define the `admissible' metrics which are rigid and also 
define a Riemannian submersion in $\calU$.  Because we are interested in precise asymptotics, it is necessary to assume
the constancy of the indicial roots associated to the Hodge Laplacians $\Delta_{k,g}$. This is equivalent to asking 
that the eigenvalues of the induced Hodge Laplacians on $(F_y,\kappa(0)|_{\phi^{-1}(y)})$ are independent of $y \in B$. 
(In practice, we only need that the eigenvalues in some fixed range $[0,C]$ are constant, but we make the stronger
assumption because it is more convenient.) 

The more general types of metrics we study are asymptotically of this form, i.e.\ $ g = g_0 + h$ where $g_0$ is admissible and $h$
is smooth (or polyhomogeneous) up to the edge and decays relative to $g_0$. Finally, amongst these asymptotically
admissible metrics we single out two special subclasses:
\begin{itemize}
\item[i)] The metric $g$ is strongly asymptotically admissible if $|h|_{g_0} = \calO(x^{b+1})$;
\item[ii)] The asymptotically admissible metric $g$ is called even if the expansion of $h$ (as described below)
contains only even powers of $x$, up to components corresponding to $x$-cross-terms, which are required 
to contain only odd powers of $x$. 
\end{itemize}
We remark that this definition of an even (asymptotically admissible) metric requires that we specify an equivalence
class of boundary defining functions $x$, where two such functions are equivalent if their quotient is a smooth
function with only even terms in its expansion up to order $b+1$. This type of even substructure arises
in many other places, for example Fefferman and Graham's ambient metric construction, as well as in
various spectral geometric and index-theoretic calculations similar to the ones here, cf. 
\cite{Albin} and \cite[\S2]{Gui}.  Similarly, strongly asymptotically admissible metrics are preserved
if we change the boundary defining function by another one which agrees with it up to order $b+1$. 
The precise definition of these `special coordinates' appears below at the end of \S 2.1. 

We may now state our main
\begin{thm}
Let $(M,g)$ be a compact simple edge space.  Then the analytic torsion, 
and hence the Ray-Singer metric $\|\cdot \|^{RS}_{(M,g)}$,
in terms of the Friedrichs extension of the Hodge Laplacian, 
is well-defined for any asymptotically admissible metric provided 
$m = \dim M$ is odd, and moreover, $\|\cdot \|^{RS}_{(M,g)}$ 
is invariant under all deformations amongst strongly asymptotically admissible metrics $g$ 
which fix the admissible rigid metric $g_0$. If $b = \dim B$ is also odd, 
then for every choice of special coordinates, $\|\cdot \|^{RS}_{(M,g)}$ is invariant under all 
deformations amongst even or strongly asymptotically admissible metrics $g$, including
those which vary the admissible rigid metric $g_0$.
\label{mainthm}
\end{thm}

The earliest study of analytic torsion and its possible relationship to combinatorial invariants in the setting of singular spaces
was by Dar \cite{Dar}. More recently, the second author studied the analytic torsion of truncated cones \cite{Vert}, with similar 
independent work by Spreafico \cite{Spr}.  Dai and Huang \cite{DH} have initiated a study of Reidemeister torsion for conic spaces,
but find that there is no obvious unique generalization of Reidemeister torsion.  The theorem above suggests that this may be
no accident, since the analytic torsion is invariant under a reasonably broad class of metrics only when $\dim B$ is odd, 
which excludes the case when $B$ is a point. We do not yet have a topological interpretation of this torsion invariant, but 
hope to return to this soon.

The main step in the proof of Theorem~\ref{mainthm} is the construction of the heat kernel for the Hodge Laplacian on $(M,g)$ 
using the methods of geometric microlocal analysis. This was carried out some time ago by Mooers \cite{Moo} for spaces with isolated 
conic singularities, and the construction here is very similar in spirit (and most details), but must be done with careful attention to 
the possible terms which can arise in the asymptotic expansion of the pointwise trace. The `geometric microlocal' method referred 
to here is the one pioneered by Melrose, see \cite{Mel-APS}, and developed by him and many others in the past few decades. It 
involves, for this problem, the careful study of the Schwartz kernel of the heat operator, as a polyhomogeneous distribution on a 
certain resolution of $M \times M \times {\mathbb R}^+$ obtained by a sequence of real blowups. Heat trace asymptotics for
simple edge spaces have also been studied by Br\"uning and Seeley \cite{BS3}; they approach these via the resolvent expansion, 
which although roughly equivalent is perhaps slightly less well adapted to the present purpose.

Our main theorem is a consequence of the following three main technical results. All notation is as above,
but we also use some terminology which will be explained later.
\begin{thm}
For each degree $k$, the (rescaled) heat kernel $H_k := e^{-t\Delta_k}$ lifts to a blown up `heat space' $M^2_h$ as a polyhomogeneous 
distribution. Its asymptotic expansions at all boundary faces of this space are determined by the indicial roots of $\Delta_k$.
In particular, it is $\rho^{-m}$ times a smooth function near the face $\td$ and $\rho^{-b-1}$ times a smooth function
near the face $\ff$, where in either case $\rho$ is a boundary defining function for the corresponding face.
\label{strhk}
\end{thm}
\begin{thm}
The trace of the heat kernel $H_k$ is a polyhomogeneous distribution on $\RR^+$ with asymptotic expansion
$$
\Tr H_k(t) \sim \sum_{\ell=0}^{\infty}A_\ell t^{\ell-\frac{m}{2}}+\sum_{\ell=0}^{\infty}C_\ell t^{\frac{\ell-b}{2}}+
\sum_{\ell\in \mathfrak{I}}G_\ell t^{\frac{\ell-b}{2}}\log t 
$$
as $t \to 0$, where
$$
\mathfrak{I}=\{\ell \in \N_0 \mid \ell + m -b \ \textup{even}\}.
$$
is the index set in the last sum. If $g$ is an even asymptotically admissible metric, then 
for every choice of special coordinates this takes the simpler form
$$
\Tr H_k(t) \sim \sum_{\ell=0}^{\infty}A_\ell t^{\ell -\frac{m}{2}}+\sum_{\ell=0}^{\infty}C_\ell t^{\ell-\frac{b}{2}}+
\sum_{\ell\in \mathfrak{I}'}G_\ell t^{\ell-\frac{b}{2}}\log t, 
$$
where $\mathfrak{I}' =\emptyset $ if $(m-b)$ is odd and $\mathfrak{I}'=\N_0 = \{0,1,2,\ldots\}$ if $(m-b)$ is even.
\label{hkasym}
\end{thm}
\begin{thm}
Suppose that $g_\mu$ is a family of asymptotically admissible metric which satisfy either conditions i) or ii) above,
with $b$ odd in the later case. Then
$$
\frac{d\,}{d\mu} \|\cdot \|^{RS}_{(M,g_{\mu})}=0.
$$
\label{metinv}
\end{thm}

Our main result, the invariance of the analytic torsion under all deformations  
amongst even asymptotically admissible edge metrics, in case of the total and 
the edge dimension being both odd, has an important 
consequence towards the Ray-Singer conjecture for a certain class of stratified 
spaces. 

More precisely, consider a smooth closed Riemannian manifold $(M^m,g)$ with a codimension two
submanifold $B^{m-2}$, and assume that $m$ (and hence $b = m-2$) is odd. View $M$ as a simple edge
space, and introduce a family of even asymptotically admissible edge metrics $g_{\A}$, where 
$\alpha \in (0,\infty)$ is the cone angle in the direction normal to $B$, with $g_{2\pi} = g$.
Since the Friedrichs extension of the Laplacian for $(M,g_{2\pi})$ coincides with the standard 
self-adjoint realization of the Laplacian on the smooth compact manifold $(M,g)$, Theorem 
\ref{mainthm} proves that for any $\alpha > 0$, 
\[
\|\cdot \|^{RS}_{(M,g_{\A})}=\|\cdot \|^{RS}_{(M,g)}.
\]
By the Cheeger-M\"uller theorem, the right side is equal to the combinatorial Reidemeister
torsion norm:
\[
\|\cdot \|^{RS}_{(M,g)}=\|\cdot \|^{\mathrm{Reid}}_M.
\]
This gives a combinatorial meaning to the analytic torsion for the edge metrics $g_\A$.
We might hope that the right hand side of this last equality agrees with the Reidemeister torsion
norm defined in terms of the intersection homology of the stratified space $(M,B)$, see \cite{Dar}.
Unfortunately, this last point is not so clear, see \cite{DH}. 

This paper is organized as follows. We begin in \S \ref{edge} with a more careful definition of the class of simple edge spaces,
the various classes of incomplete edge metrics mentioned above, and an examination of the structure of the Hodge Laplacian
for these metrics. We also consider the asymptotics of solutions to $\Delta w = 0$  where $w$ is in the maximal extension 
of this operator, and use this to characterize the Friedrichs extension. \S \ref{construct} contains the parametrix construction
for the heat kernel as an element of the calculus of heat operators on $M$; we also introduce the even subcalculus and
show that it contains $H_k$ when $g$ is even. These results lead directly to a description of the asymptotics of the
heat trace in \S 4; the proof that the analytic torsion is well-defined and the computation of its variation is given in
\S \ref{torsion}. The Appendix contains the proof of the composition formula for the heat calculus, as well as the even
subcalculus. 

\bigskip

\section*{Acknowledgements}

The authors would like to thank Pierre Albin, Matthias Lesch and Werner M\"uller for helpful discussions and for
their interest in this work. 
We owe particular thanks to the anonymous referee who read the manuscript exceptionally carefully and made many useful
remarks and suggestions. The first author is also grateful to Xianzhe Dai for many illuminating conversations about analytic torsion, 
while the second author wishes to thank Daniel Grieser and Andras Vasy for introducing him to elements of microlocal analysis, 
as well as to Stanford University for its hospitality. He gratefully acknowledges financial support by the German Research 
Foundation (DFG) during his postdoctoral studies at Stanford University and the Hausdorff Center for Mathematics in Bonn.  
The first author was also partially supported by the NSF grant DMS-0805529 and DMS-1105050, and also by the CNRS at 
the Universit\'e de Nantes, Laboratoire Jean Leray.

\section{Simple edge spaces and the Hodge Laplacian}\label{edge}

We begin by introducing the class of stratified spaces with simple edges and incomplete edge metrics, and then
describe the structure of the Hodge Laplace operator on such spaces near the edge. We then define the class of
polyhomogeneous distributions, and finally describe the Friedrichs extension of the Hodge Laplacian.

\subsection{Simple edge spaces}
Let $\overline{M}$ be a compact stratified pseudomanifold, with top-dimensional (open, dense) stratum $M$, an
open $m$-manifold, and a single lower-dimensional stratum $B$. There are a set of axioms, reviewed in \cite{ALMP}, which
regulate the differential topological structure of a stratified space. For our purposes, there are two main consequences of these
axioms. First, in our setting there are no higher depth strata on the frontier of $B$, so it is a compact smooth manifold; 
we could easily treat the 
case where $B$ has components of varying dimension, but for simplicity we assume that $B$ is connected, with $\dim B = b$. 
Second, there are a neighbourhood $\calU$ of $B$ in $\overline{M}$, a `radial' function $x$ defined on $M \cap \calU$ and a smooth 
projection $\phi: \calU \to B$, which is a submersion on $\calU \cap M$, such that the preimages $\phi^{-1}(q)$ are all diffeomorphic to
truncated open cones $C(F)$ over a compact smooth manifold $F$, with $\dim F = f$, and where the restriction of the function $x$ to
each fibre $\phi^{-1}(q)$ is a radial function on that cone.
We suppose that the level set $\{x=1\}$ corresponds to the `outer'
boundary of the neighbourhood $\calU$, and we denote this manifold by $Y$. Thus $Y$ is smooth and compact, and is the
total space of a fibration $\phi: Y \to B$ with fibre $F$. Figure 1 illustrates this cone bundle structure.
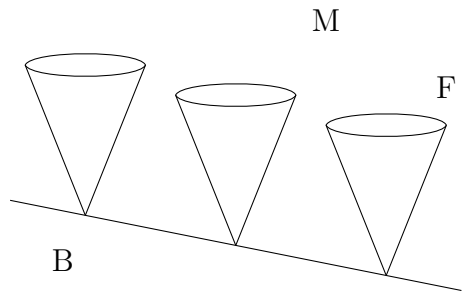
\begin{figure}[h]
\begin{center}
\begin{tikzpicture}
\draw (-3,0.6) -- (3,-0.6);

\draw (-2.8,2.4) -- (-2,0.4);
\draw (-1.2,2.4) -- (-2,0.4);
\draw (-2.8,2.4) .. controls (-2.7,2.6) and (-1.3,2.6) .. (-1.2,2.4);
\draw (-2.8,2.4) .. controls (-2.7,2.2) and (-1.3,2.2) .. (-1.2,2.4);

\draw (-0.8,2) -- (0,0);
\draw (0.8,2) -- (0,0);
\draw (-0.8,2) .. controls (-0.7,2.2) and (0.7,2.2) .. (0.8,2);
\draw (-0.8,2) .. controls (-0.7,1.8) and (0.7,1.8) .. (0.8,2);

\draw (2.8,1.6) -- (2,-0.4);
\draw (1.2,1.6) -- (2,-0.4);
\draw (2.8,1.6) .. controls (2.7,1.8) and (1.3,1.8) .. (1.2,1.6);
\draw (2.8,1.6) .. controls (2.7,1.4) and (1.3,1.4) .. (1.2,1.6);

\node at (-2.3,-0.2) {\large{B}};
\node at (2.8,2.1) {\large{F}};
\node at (1.2,3) {\large{M}};

\end{tikzpicture}
\end{center}
\label{edge-picture}
\caption{The tubular neighbourhood $\calU$ in $M$.}
\end{figure} 

The class of metrics $g$ we consider on such a space are the ones which restrict, on each conical fibre, to be asymptotically conic.
To be more precise, let us say that $g$ is {\it rigid} if it is arbitrary away from the open set $\calU$, but that in $\calU$ it has 
the following form: there is a smooth Riemannian metric $\ha$ on $B$ and a symmetric $2$-tensor $\ka$ on $Y$ which
restricts to a metric on each fibre $F$ such that 
\[
\left. g \right|_{\calU} = dx^2 + \phi^*\ha +  x^2 \ka.
\]
Hence in $\calU$, the induced metric on each fibre $\phi^{-1}(q)$ is an exact warped-product conic metric. In fact, it is necessary
for us to work with a slightly more restricted class of metrics: a metric $g$ is called {\it admissible} 
if it is rigid and in addition $\phi: (Y, \left. g \right|_Y) \to (B,\ha)$ is a Riemannian submersion. Recall
that this means the following: if $p \in Y$, then $T_pY$ splits into vertical and horizontal subspaces, $T^V_p Y \oplus
T^H_p Y$, where by definition $T^V_pY$ is the tangent space to the fibre of $\phi$ through $p$ and $T^H_p Y$ is the
orthogonal complement of this subspace. The new condition is that the restriction of the tensor $\ka$ to $T^H_p Y$ vanishes. 

More generally, the metric $g$ is asymptotically admissible if $g = g_0 + h$ where $g_0$ is admissible and
$|h|_{g_0} \to 0$ at $B$. We shall typically assume that $h$ has a polyhomogeneous expansion in powers of $x$
(see \S 2.4), or even more strongly that it is smooth to $x=0$. In order to describe this accurately, 
it is convenient to pass from the singular space $\overline{M}$ to its resolution $\widetilde{M}$, the manifold
with boundary obtained by `blowing up' the stratum $B$. In concrete terms, this corresponds to replacing each
conic fibre $C(F) = [0,1)_x \times F/ \sim$ (where $(0,q) \sim (0,q')$ for any points $q, q' \in F$) with the
cylinder $[0,1)_x \times F$, and correspondingly, replacing the cone bundle neighbourhood $\calU$ with the 
associated bundle of cylinders. The $\calC^\infty$ structure on $\wt{M}$ is the natural one induced from the
cylindrical fibres. Notice that $\partial \widetilde{M} = Y$. This resolution 
process is described in greater detail in \cite{Maz-edge}; see also \cite{ALMP}. 

The advantage of considering objects on $\wt{M}$ rather than $\overline{M}$ is that there is now a clear
meaning for a function or tensor to be smooth to $x=0$. Thus we now consider metrics
\[
\left. g \right|_{\calU} = dx^2 + \phi^*\ha(x) +  x^2 \ka(x) + \ \textup{higher orders}.
\]
where $\ha(x)$ and $\ka(x)$ depend smoothly on $x \in [0,1)$. We  can always decompose any such $g$ 
as $g_0 + h$ where $g_0 = dx^2 + \phi^* \ha(0) + x^2 \ka(0)$. Then $g$ is called strongly asymptotically
admissible if it is smooth and if $|h|_{g_0} = \calO(x^{1+b})$. It is called even if $\ha(x)$ and $\ka(x)$
are smooth as functions of $x^2$ rather than just $x$, and if in addition $h$ has a smooth expansion with
only even powers of $x$, up to components corresponding to $x$-cross-terms, which are required to contain only 
odd powers of $x$. 

Evenness of $g$ depends on restricting to a particular `even' equivalence class of coordinate charts near the edge. 
If $g$ is even with respect to one coordinate system $(x,y,z)$, then another coordinate system $(\wx, \wy, \wz)$ 
is in this equivalence class if $\wx/x$, $\wy$ and $\wz$ are all even functions of $x$, with coefficients in the 
expansions depending smoothly on $y$ and $z$. We assume henceforth that this restriction is made, often
without comment.  We noted already in the introduction that metric independence of analytic torsion requires
only that metric and coordinate charts be even up to the order $b+1$.

\subsection{Edge operators}
The paper \cite{Maz-edge} discusses in great detail the theory of elliptic operators on simple edge spaces. We review
some of the relevant material in this theory.

Consider local coordinates $(x,y,z)$ on $\widetilde{M}$ near the boundary face $Y$, where $x$ is the radial coordinate, $y$ is
the lift of a local coordinate system on $B$ and $z$ restricts to coordinates on each fibre $F$. Now define the class of edge
vector fields $\calV_e$ on $\widetilde{M}$: these are the vector fields on this space which are smooth even at $x=0$ and which are
tangent to the fibres of $Y$ at this boundary face. In this local coordinate system, any element of $\calV_e$ can be written as
a sum of smooth multiples of the basic generators $x\del_x$, $x\del_{y_i}$ and $\del_{z_j}$, which we write as
\[
\calV_e = \mbox{Span}_{\calC^\infty}\, \{x\partial_x, x\partial_{y_1}, ..., x\partial_{y_b}, \partial_{z_1},..., \partial_{z_f}\}.
\]
These are the basic objects of the edge theory. Notice that if $g$ is a smooth asymptotically admissible metric on $M$, and 
if $V, W \in \calV_e$, then $x^{-2}g(V,W) \in \calC^\infty$; in particular, the length of any edge vector field decays at least 
like $x$ as $x \to 0$. 

We next introduce the class of differential edge operators $\textup{Diff}^*_e(M)$. By definition, $L \in \textup{Diff}^*_e(M)$
if it can be written locally as a sum of products of elements of $\calV_e$, with coefficients in $\calC^\infty(\wt{M})$. Thus
\[
L=\sum_{j+|\A|+|\beta|\leq m} a_{j,\A,\beta}(x,y,z)(x\partial_x)^j(x\partial_y)^{\A}\partial_z^{\beta},
\]
with each $a_{j,\alpha,\beta}$ smooth up to $x=0$. More generally, if $L$ acts between sections of two vector bundles, then $L$ 
has this form with respect to suitable local trivializations, where each $a_{j,\alpha,\beta}$ is matrix-valued. 
The operator is called edge elliptic if its edge symbol
\[
{}^e \sigma_m(L)(x,y,z;\xi,\eta,\zeta) :=  \sum_{j + |\A| + |\be| = m} a_{j,\A,\be}(x,y,z) \xi^j \eta^{\A} \zeta^{\be}
\]
is nonvanishing (or invertible, if matrix-valued), for $(\xi,\eta,\zeta) \neq (0,0,0)$. This has an invariant meaning as a
function on the so-called edge cotangent bundle ${}^e T^* \widetilde{M}$ which is homogeneous of degree $m$ on the fibres.

There is an entire zoology of objects associated to edge geometry, but in the interests of space, we refer to the papers cited
above for more on all of this. 

The fundamental tool in the analysis of elliptic edge operators is the space of pseudodifferential edge operators $\Psi^*_e(\wt{M})$. 
These operators are standard classical pseudodifferential operators in the interior, but are adapted to the degeneracy structure
of elements of $\textup{Diff}^*_e$ and form a suitably broad class of operators so that, in many cases, an elliptic differential
edge operator $L$ has a pseudodifferential edge parametrix $G$ such that both $GL - I$ and $LG - I$ are compact on certain
natural function spaces. This is the main content of \cite{Maz-edge}. 

One key idea in this theory is the use of two different model operators, which provide fundamental information about the
edge elliptic operator $L$ in the parametrix construction. The first of these is the indicial operator, 
\[
I(L) = \sum_{j + |\be| \leq m} a_{j,0,\be}(0,y_0,z) (s\del_s)^j \del_z^\be,
\]
acting on functions on $\RR^+ \times F$, where $y_0$ is some fixed point on $B$; by taking the Mellin
transform in $s$, this reduces to the indicial family
\[
I_\zeta(L) = \sum_{j + |\be| \leq m} a_{j,0,\be} (0,y_0,z) \zeta^j \del_z^\be,
\]
which is a holomorphic family of unbounded Fredholm operators on $L^2(F)$. Values of $\zeta$ for which 
$I_\zeta(L)$ is not invertible are called indicial roots of $L$; the indicial family is a generalization of the resolvent
(of, say, $I_0(L)$), and the indicial roots then play the role of the eigenvalues of this operator.
Using the ellipticity of $L$, it may be shown that the set of indicial roots is discrete in ${\mathbb C}$ and that the 
corresponding solutions of $I_\zeta(L) \phi(z)=0$ lie in $\calC^\infty(F)$. The second model is the normal operator, 
\[
N(L) = \sum_{j + |\A| + |\be| \leq m} a_{j,\A,\be}(0,y_0,z) (s\del_s)^j (s\del_u)^{\A} \del_z^{\be}.
\]
Here $(s,u) \in {\mathbb R}^+ \times {\mathbb R}^b$ are linear variables on a half-space, which should be thought of 
as the inward pointing normal space to the fibre of $Y$ through $(0,y_0,z)$. This seems to have almost the same complexity
as $L$, but since it is translation invariant in $u$ and dilation invariant in $(s,u)$ jointly, it can be reduced by Fourier
transform and rescaling to an operator on $\RR^+ \times F$ which is only slightly more complicated than $I(L)$ 
to analyze. 

One main result in the theory is that if $N(L)$ is invertible (on some fixed weighted $L^2$ space) for each $y_0 \in B$, then $L$ 
itself is Fredholm on the corresponding weighted $L^2$ space on $M$. The indicial roots of $L$ determine a discrete set of
weights such that $L$ does not have closed range when acting on the corresponding weighted $L^2$ space.  In any case,
the point is simply that the inverse of $N(L)$ is the main ingredient in the construction of a 
parametrix for $L$. Something similar is true for the construction of a heat kernel parametrix for $L$. 

We conclude this general discussion by observing that although the Laplacian and other natural elliptic operators on a simple 
edge space with incomplete edge metric are not quite edge operators, they are very closely related and can be studied by this 
edge theory. More specifically, if $g$ is in one of the classes of incomplete edge metrics above and $\Delta_g$ its scalar Laplacian,
then $\Delta_g = x^{-2}L$ where $L \in \textup{Diff}^{\,2}_e(M)$. The analogous statement is true for the Hodge Laplacian 
provided we trivialize the form bundles appropriately. We have already alluded to the edge cotangent bundle ${}^e T^*M$.
In the local coordinates $(x,y,z)$, this bundle has a local smooth basis of sections
\[
\left\{\frac{dx}{x}, \frac{dy^1}{x}, \ldots, \frac{dy^b}{x}, dz^1, \ldots, dz^f\right\}.
\]
The edge $k$-form bundle ${}^e \Lambda^k (M)$ is simply the $k^{\mathrm{th}}$ exterior power of this bundle, hence
is generated locally by $k$-fold wedge products of these sections. The correct assertion, then, is that
$\Delta_{g,k}$, as an operator acting on sections of ${}^e \Lambda^k(M)$, has the form $x^{-2}L$ where $L \in \textup{Diff}^2_e(M;
{}^e\Lambda^k(M))$. We omit the proofs of these facts; they can either be checked by direct computation or else 
inferred using the `naturality' of edge structures, cf.\ \cite{Maz-edge}. 

Similarly, we define ${}^{ie} T^*M$ to be spanned by a local smooth basis of sections
\[
\left\{ dx, dy^1, \ldots, dy^b, x dz^1, \ldots, x dz^f\right\}.
\]
By definition, the incomplete edge $k$-form bundle ${}^{ie} \Lambda^k (M)$ is then the 
$k^{\mathrm{th}}$ exterior power of ${}^{ie} T^*M$. Denote smooth sections of
${}^{ie} \Lambda^k (M)$ by ${}^{ie} \Omega^k (M)$.

\subsection{The Hodge Laplacian}\label{rescaling-section}
We now consider the structure of the Hodge Laplacian of an incomplete edge metric.  More specifically, we describe some part
of the structure of the normal operator $N(x^2\Delta_p)$ and of its indicial roots. The material here is drawn from
\cite{HM}, which in turn summarizes and presents in unified fashion various results proved by Bismut-Cheeger and Dai.

The normal operator for $x^2 \Delta_p$ at any $y_0 \in B$ acts on $p$-forms on the model edge $\RR^b \times C(F) = 
\RR^+_s \times \RR^b_u \times F$ with incomplete edge metric $\overline{g} = ds^2 + s^2 \ka(0) + |du|^2$, and is naturally 
identified with $s^2$ times the Hodge Laplacian for that model metric. 

The first step is to consider the structure of the operator induced on the hypersurface $S = \{s=1\} = \RR^b \times F$. 
For this, note that $T(\RR^b \times F)$ splits into the sum of a `vertical' and `horizontal' subspace, where the first is the 
tangent space to the $F$ factor and the second is the tangent space to the Euclidean factor. This splitting is orthogonal, and
induces a bigrading 
\[
\Lambda^p(S) = \bigoplus_{j + \ell = p} \Lambda^j(\RR^b) \otimes \Lambda^\ell (F) := \bigoplus_{j+\ell = p}  \Lambda^{j,\ell}(S). 
\]
Let ${}^e \Omega^{j,\ell}(S)$ denote the space of sections of the corresponding summand in this bundle decomposition. 
The differentials and codifferentials on these factors satisfy
\begin{align*}
& d_{\RR^b}: \Omega^{j,\ell}(S) \to \Omega^{j+1,\ell}(S),  \qquad \qquad & d_F: \Omega^{j,\ell}(S) \to  \Omega^{j,\ell+1}(S) \\
& \delta_{\RR^b}: \Omega^{j,\ell}(S) \to \Omega^{j-1,\ell}(S),  \qquad \qquad & \delta_F: \Omega^{j,\ell}(S) \to 
\Omega^{j,\ell-1}(S).
\end{align*}
Furthermore, $d_S = d_{\RR^b} + d_F$ and $\delta_S = \delta_{\RR^b} + \delta_F$. It is not hard to check that the 
differential and codifferential induced on each level set $S_a = \{s=a\}$ has the form 
\[
d_{S_a} = d_{\RR^b} + d_F, \qquad \delta_{S_a} = \delta_{\RR^b} + a^{-2} \delta_F.
\]

This can all be assembled into an expression for $N(x^2\Delta_p)$ on $\RR^b \times C(F)$.  In order to simplify various calculations
below, we write this operator using a rescaling of the form bundles, employed also by Br\"uning-Seeley \cite{BS}
and in slightly different form in \cite{HM}.  Thus, for each $j,\ell$ with $j + \ell = p$, define
\begin{align*}
& \phi_{j,\ell}: \calC^{\infty}(\RR^+, \Omega^{j,\ell-1}(S) \oplus \Omega^{j, \ell}(S)) \rightarrow {}^{ie}\Omega^p(\RR^b \times C(F)), \\
& \qquad (\eta, \mu) \longmapsto s^{\ell-1-f/2}\eta \wedge ds + s^{\ell-f/2}\mu, 
\end{align*}
and denote by $\Phi_p$ the sum of these maps over all $j+\ell = p$. It is not hard to check that
\begin{equation}
\label{unitary}
\begin{split}
\Phi_p: L^2\left(\RR^+, L^2\big(\bigoplus_{j+\ell=p} \Omega^{j,\ell-1}(S) \oplus \Omega^{j,\ell}(S), \kappa(0) + |du|^2 \big) , ds\right)&
\\ \longrightarrow L^2(\,{}^{ie}\Omega^p(\RR^b \times C(F)), \overline{g}),&
\end{split}
\end{equation}
is an isometry, and a calculation yields
\begin{align}\label{laplace}
\Phi_p^{-1}\circ \left[s^{-2}N(x^2\Delta_p)\right]\circ \Phi_p=\left(-\frac{\del^2}{\del s^2}+\frac{1}{s^2}(A-1/4)\right) + \Delta_B, 
\end{align}
where $A$ is the nonnegative self-adjoint operator on $\Lambda^{\ell-1}(F) \oplus \Lambda^{\ell}(F)$ given by 
\begin{align}\label{a}
A=\left(\begin{array}{cc}\Delta_{\ell-1,F} + (\ell-(f +3)/2)^2 & 2(-1)^{\ell}\, \delta_{\ell,F}\\ 2(-1)^{\ell}\, d_{\ell-1,F}& \Delta_{\ell,F}+ 
(\ell-(f +1)/2)^2\end{array}\right).
\end{align}
If $(M,g)$ is a simple edge space with an asymptotically admissible edge metric, then we can define a similar
rescaling $\Phi$ using powers of the defining function $x$ locally in the neighbourhood $\calU$ of $B$. Assuming that $x$ is
smooth on $M \setminus B$ and equals $1$ away from $\calU$, then this rescaling can be extended trivially
to the rest of $M$. Rescalings on different local coordinate neighborhoods are equivalent up to a diffeomorphism. 
Conjugating by $\Phi$ as before, then this rescaled operator $\Delta_p^\Phi$ is a perturbation of 
\eqref{laplace} with higher order correction terms determined by the curvature of the Riemannian submersion $\phi: Y \to B$ 
and the second fundamental forms of the fibres $F$ in $Y$. We write this rescaled operator simply as 
$\Delta_p$ if there is no danger of confusion.

One reason for this transformation is that the indicial roots of $\Delta_p^\Phi$ have a particularly simple form:  
writing the eigenvalues of $A$ as $\nu_j^2$, with corresponding eigenform $\phi_j$, the
corresponding indicial roots of $\Delta_p$ are the roots of the quadratic equation
\begin{equation}
-\gamma (\gamma - 1) + \nu_j^2 - \frac{1}{4} = 0 \Leftrightarrow \gamma_j^{\pm} = \frac{1}{2} \pm \frac{1}{2}
\sqrt{1 - 1 + 4\nu_j^2} = \frac12 \pm \nu_j, 
\label{indroots}
\end{equation}
where we make the convention that every $\nu_j \geq 0$. If $\nu_0 = 0$, then the corresponding solutions are $x^{1/2}\phi_0$
and $x^{1/2}(\log x) \phi_0$. Note that $s^{\frac12 - \nu_j}\phi_j \notin L^2$ when $\nu_j \geq 1$.

\subsection{Polyhomogeneity}
Before continuing, we must introduce the spaces of conormal and polyhomogeneous distributions on a manifold with 
corners $W$. These are generalizations of and replacements for $\calC^\infty$ functions in this setting, and are needed
simply because solutions of edge elliptic and parabolic equations tend to have this form.  Briefly, a function is 
conormal if it has stable regularity with respect to differentiations by arbitrary smooth vector fields which are 
tangent to all boundaries and corners of $W$; it is polyhomogeneous if it has an asymptotic expansion 
at all boundary faces, and a product type expansion at all corners, in terms of powers of the boundary defining
functions for the hypersurface faces, with all coefficients depending smoothly on the tangential variables.  
The exponents which appear in these expansions may be arbitrary complex numbers, and we also allow 
nonnegative integer powers of the logarithms of the defining functions as factors.   

To state all of this more formally, let $W$ be a compact manifold with corners, with all boundary faces embedded, and 
$\{H_i\}_{i=1}^N$ the set of all hypersurface boundaries of $W$. As part of our definition of manifold with corners,
we require that each $H_i$ be embedded, so it has a boundary defining function $\rho_i$, i.e.\ $\rho_i$ is a smooth
nonnegative function on  $W$ which vanishes simply on $H_i$ and is strictly positive on $W \setminus H_i$. 
For any multi-index $\lambda= (\lambda_1,\ldots, \lambda_N) \in \RR^N$ we write $\rho^\lambda = 
\rho_1^{\lambda_1} \ldots \rho_N^{\lambda_N}$.  Finally, denote by $\calV_b(W)$ the space of all smooth vector fields on 
$W$ which are unconstrained in the interior but which lie tangent to all boundary faces. Thus, if $p$ is a point
on a codimension $k$ corner of $W$, then there exists a local coordinate chart near $p$, $(x_1, \ldots, x_k, y_1,
\ldots, y_{n-k})$ with each $x_i \in [0,\e)$ and $y_j \in (-\e, \e)$, and very similarly to our definition of $\calV_e$ in \S 2.2,
\[
\calV_b(W) = \mbox{Span}_{\calC^\infty}\, \{ x_i \del_{x_j}, \del_{y_s},\ i,j = 1, \ldots, k,\ s = 1, \ldots, n-k \}.
\]

\begin{defn}\label{phg}
A distribution $w$ on $W$ is said to be conormal of order $\lambda$ if it is of stable regularity with respect to
elements of $\calV_b(W)$, i.e.\ if it lies in the space
\[
\calA^\lambda(W) = \{ w \in \rho^\lambda L^\infty(W):  V_1 \ldots V_\ell w \in \rho^\lambda L^\infty(W)
\ \forall\, V_j \in \calV_b(W), \ell \geq 0\}.
\]
We also write $\calA^*(W) = \cup_\lambda \calA^\lambda(W)$. 
Note that any conormal distribution is smooth in the interior of $W$.  Furthermore, $L^\infty$ could
be replaced by any other fixed space, e.g.\ some (polynomially) weighted $L^2$ or $L^p$ space; the union of 
these spaces over all $\lambda$ is the same in any case, only the weight and regularity scale would change.

\medskip

Next, an index set $E = \{(\gamma,p)\} \subset {\mathbb C} \times {\mathbb N}$ is the set of exponents 
associated to an expansion at the face $H_i$. We require that it satisfies the following hypotheses:
\begin{enumerate}
\item Each half-plane $\textup{Re} \, \zeta < C$ contains only finitely many $\gamma$;
\item For each $\gamma$, there is a $P(\gamma) \in {\mathbb N}_0$ such that $(\gamma,p) \in E_i$ iff 
$0 \leq p \leq P(\gamma) < \infty$;
\item If $(\gamma,p) \in E_i$, then $(\gamma+j,p) \in E_i$ for all $j \in {\mathbb N}$.
\end{enumerate}
An index family $\calE = (E_1, \ldots, E_N)$ is an $N$-tuple of index sets, where $E_j$ is associated to the
face $H_j$ of $W$.

\medskip

Finally, we say that $w$ is polyhomogeneous on $W$ with index family $\calE$ if it has expansions at 
the various boundary hypersurfaces with exponents determined by the index family $\calE$:
\[
\calA_{\phg}^\calE(W) = \{w\in \calA^*(W): \ w \sim \sum_{(\gamma,p) \in E_i} a_{\gamma,p}^{(i)}
\rho_i^{\gamma} (\log \rho_i)^p \ \ \mbox{near $H_i$}\},
\]
where all coefficients are themselves polyhomogeneous, $a_{\gamma,p}^{(i)} \in \calA_{\phg}^{\calE^{(i)}}(H_i)$.
The index family $\calE^{(i)}$ is the obvious one induced from $\calE$ at the boundary faces 
$H_i \cap H_j$ of $H_i$. 
\end{defn}

It is a simple consequence of this definition that if $w$ is polyhomogeneous, then it has a product type
expansion at any corner $H_{i_1} \cap \ldots \cap H_{i_\ell}$ of $W$ of the form
\[
w \sim \sum a_{\gamma,p} \rho^\gamma (\log \rho)^p,
\]
where now $\gamma$ and $p$ are multi-indices, and with coefficient functions conormal on that corner. 
We recall also that by saying that the expansion for $w$ is asymptotic, we mean that the difference between 
$w$ and any finite portion of the expansion vanishes at the rate of the next term in the expansion, with a 
corresponding property for all higher derivatives.

In the next subsection we shall also encounter distributions which are conormal and `partially polyhomogeneous',
i.e.\ they are conormal and have a finite expansion up to some order of decay, with a remainder term which
is only conormal. We do not introduce special notation for these spaces.
 
We refer to \cite{Maz-edge} for more details about polyhomogeneity. 

\subsection{The Friedrichs extension}\label{bc}
Since the compact simple edge space $(M,g)$ is incomplete, the Hodge Laplacian may not be essentially self-adjoint
on the core domain $\calC^\infty_0\Omega^p(M)$, so we must consider how to impose boundary conditions at the edge to
obtain closed, or even better, self-adjoint, extensions.  For spaces with isolated conic singularities, this was first
accomplished by Cheeger \cite{Ch2}.  Further and more systematic studies in the conic setting for general Dirac-type operators
appear in \cite{Lesch} and \cite{GM}, and see \cite{Moo}, \cite{KLP} for results about the associated heat equation. 
That setting is tractable because the extension problem is finite dimensional. When the edge has positive dimension, 
the analysis needed to carry this out is more intricate. We restrict attention to the Friedrichs extension since this
requires much less machinery to define and characterize. Other self-adjoint extensions have been studied by the 
second author, with Bahuaud and Dryden in \cite{BDV}.

Let $\Delta_p^\Phi = \Delta_p$ denote the rescaled Hodge Laplace operator acting on differential forms of degree $p$
on the compact simple edge space $M$.  Consider the space of $L^2$ forms $L^2\Omega^p(M)$, 
with respect to any choice of (polyhomogeneous) incomplete iterated edge metric on $M$, as well as the associated 
edge Sobolev spaces
\[
H^\ell_e \Omega^p(M) :=\{u\in L^2\Omega^p(M) \mid V_1\cdots V_j u \in L^2\Omega^p(M) \ \mbox{for}\ V_i \in \calV_e\ 
\mbox{and for any}\  j\leq \ell\}.
\]
We often use $H^\infty_e$  to denote the intersection over all $\ell$ of the spaces $H^\ell_e$. 
The maximal domain of $\Delta_p$ is, by definition,
\[
\calD_{\max}(\Delta_p) := \{ u\in L^2\Omega^p(M) \mid  \Delta_p u \in L^2\Omega^p(M) \}, 
\] 
where $\Delta_p u\in L^2$ is initially understood in the distributional sense. Similarly, the minimal domain of $\Delta_p$ 
is defined as
\begin{gather*}
\calD_{\min}(\Delta_p) := \{ u \in \calD_{\max}(\Delta_p)  \mid \exists\,  u_j \in \calC^\infty_0\Omega^p\ \mbox{such that}  \\ 
u_j \to u\ \mbox{and}\ \Delta_p u_j \to \Delta_p u\ \mbox{both in}\ L^2\Omega^p \}.
\end{gather*}
These are the domains of the maximal and minimal extensions of $\Delta_p$ on the core domain $\calC^\infty_0\Omega^p(M)$. 
The set of all closed extensions of $\Delta_p$ is in bijective correspondence with the closed subspaces of the quotient 
$\calD_{\max}/\calD_{\min}$; furthermore, since $\Delta_p$ is symmetric on the core domain, self-adjoint extensions are in 
bijective correspondence with the subspaces of this quotient which are Lagrangian with respect to a certain natural symplectic 
form induced from the boundary contributions in an integration by parts formula, see \cite{Lesch} and \cite{GM}. 

This motivates the problem of characterizing elements in the maximal domain. For spaces with isolated conic singularities, this is
straightforward and leads to an explicit parametrization of all closed and self-adjoint extensions of $\Delta_p$.  
For spaces with simple edges singularities, however, the complete characterization remains functional analytic, 
although there are many explicit extensions which parallel the definitions in the conic setting.
Different choices of closed extensions in either case correspond to what are sometimes called ideal boundary 
conditions at the singular stratum. We begin with a result which holds in both the conic and edge settings.

\begin{lemma}\label{max}
Let $(M,g)$ be compact with simple edge singularity. Then 
\begin{align*}
\calD_{\max}(\Delta_p)&\subset H^\infty_e\Omega^p(M) + x^2 H^2_e\Omega^p(M) \subset H^2_e\Omega^p(M),  \\
\calD_{\min}(\Delta_p) &\subset x^{2-\epsilon} H^2_e\Omega^p(M)\ \mbox{for any} \ \epsilon >0.
\end{align*}
If $M$ has only conic singularities, then any $w \in \calD_{\max}(\Delta_p)$ admits an asymptotic expansion
\begin{equation}
w \, \sim \, \tilde{w} + \sum_{j=1}^N x^{\frac12 + \nu_j}  a^+_j(z) + \left\{
\begin{aligned} &x^{\frac12 - \nu_j} a^-_j(z), \quad \, \, \, & \nu_j\neq 0 \\ & x^{\frac12} \log x \, a^-_j(z), \, & \nu_j=0 
\end{aligned} \right.  
\qquad \tilde{w} \in \calD_{\min}(\Delta_p),
\end{equation}
where the numbers $1/2 \pm \nu_j$ are the indicial roots of $\Delta_p$ with $\nu_j\in [0,1)$ and 
each $a_j^\pm(z)\in \calC^\infty(F)$ is a solution of the corresponding indicial operator. If $M$ has a simple edge 
singularity, then $w$ admits an analogous expansion with coefficients $a_j^\pm(y,z)$ now depending on both $y$ and 
$z$, but this is an asymptotic expansion only in a weak sense, i.e.\ there is an expansion of the pairing
$\int_B w(x,y,z) \chi(y)\, dy$ for any test functions on $B$. 
\end{lemma}
\begin{remark}
We refer to  \cite{Maz-edge} for a careful explanation of such weak expansions. It is the failure of the expansion to hold 
in a strong sense which makes this result more difficult to use for simple edge spaces with $\dim B > 0$.  See also 
\cite{HM} and \cite{ALMP} for the case where there are no indicial roots in the critical range, so that $\calD_{\max} = \calD_{\min}$. 
\end{remark}
\begin{proof}
Consider any $w\in L^2\Omega^p(M)$ with $\Delta_p w=f\in L^2\Omega^p(M)$. Then $x^2 \Delta_p w = x^2 f \in 
x^2 L^2\Omega^p(M)$, which is useful since the operator $x^2 \Delta_p$ is an elliptic differential edge operator of order 
$2$, which we denote by $\calL$.  Applying \cite[Theorem 3.8]{Maz-edge}, we obtain a parametrix 
$\calB \in \Psi^{-2}_e(M)$ for $\calL$ in the small edge calculus, which satisfies
\[
\calB \calL =I-\calR, \ \calR \in \Psi^{-\infty}_e(M).
\]
Corollary 3.23 in \cite{Maz-edge} asserts 
\begin{align*}
\calB: x^s H^\ell\Omega^p(M) \to x^s H^{\ell+2}\Omega^p(M), \\
\calR: x^s H^\ell\Omega^p(M) \to x^s H^\infty_e\Omega^p(M),
\end{align*}
for any $s \in \RR$, $\ell \in {\mathbb N}$. Applying $\calB$ to $\calL w = x^2 f$ yields $w= \calB (x^2 f) + \calR (w)$, 
and this lies in $H^\infty_e + x^2 H^2_e$ as claimed. 

If $M$ has only isolated conic singularities, then there is a considerable sharpening of the parametrix construction,
see \cite{Maz-edge} or \cite{Mel-APS}. It is then possible to choose a much better parametrix $\calB$ which has
the property that the remainder term $\calR = I - \calB \calL$ maps an arbitrary element of $L^2\Omega^p(M)$ into a
polyhomogeneous distribution. We do not describe that construction, but it is very similar in spirit (and simpler
than) the construction for the heat parametrix in \S 4 below.  In any case, granting this, then the corresponding
equation $u = \calB (x^2 f) + \calR w$ now gives that $u \in x^2 H^2_b(M) + \calA_{\phg}$, i.e.\ it has a partial expansion 
up to order $x^2$ as in the statement of the lemma.  Once we \emph{know} that $u$ has such a partial expansion,
we can determine the structure of the terms which vanish less quickly than $x^2$ simply by substituting this
expansion into the equation and calculating formally. The only possible terms $x^\gamma \phi(z)$, $\gamma < 2$,
for which $\calL (x^\gamma \phi(z)) = \calO(x^2)$ are those with $\gamma$ indicial and $\phi$ a corresponding
solution. 

The proof that such a partial expansion holds in the weak sense in the simple edge case uses the Mellin transform,
and is explained carefully in \S 7 of \cite{Maz-edge}.

If $w \in \calD_{\max}$, then a straightforward calculation shows that any function which can be approximated
in the graph norm cannot have any terms of the form $x^\gamma a(y,z)$ (even if just in a weak expansion).
This proves the claim about the minimal domain.
\end{proof}

To characterize the Friedrichs extension of the Hodge Laplacian, observe that $\Delta = \bigoplus_p \Delta_p$
factors as $\Delta=D\circ D$, where $D = d + \delta$ is the Gau{\ss}-Bonnet operator on $\Omega^*(M)$.  
Letting $D$ act on sections of ${}^e\Omega^*(M)$, as above, then $xD$ is an elliptic edge operator of order $1$. 
Its normal operator $N(xD)_{y_0}$ at any $y_0 \in Y$ acts on the model edge $\RR^b \times C(F)$ and in fact $s^{-1}N(xD)_{y_0}$ is naturally
identified with the Gau\ss-Bonnet operator $D_{\overline{g}}$ for the model incomplete edge metric 
$\overline{g} = ds^2 + s^2 \ka + |du|^2$. Conjugating by the unitary transformation $\Phi$ from
\eqref{unitary} we calculate that, according to the splitting $\Omega^* = \Omega^{\mathrm{even}} \oplus
\Omega^{\mathrm{odd}}$, 
\begin{equation}\label{D1}
D_{\overline{g}}= 
\begin{pmatrix}
0  & \left(- \frac{\partial}{\partial s}+\frac{1}{s}P\right)+D_{\RR^b}, \\ 
\left(\frac{\partial}{\partial s}+\frac{1}{s}P\right)+D_{\RR^b}, & 0
\end{pmatrix}
\end{equation}
and so
\begin{equation} \label{D2}
D_{\overline{g}}^2 = 
\begin{pmatrix}
-\frac{\del^2}{\del s^2}+\frac{1}{s^2}\left(A_+ -1/4\right) + \Delta_{\RR^b} & 0 \\
0 & -\frac{\del^2}{\del s^2}+\frac{1}{s^2}\left(A_- -1/4\right) + \Delta_{\RR^b}
\end{pmatrix}.
\end{equation}
Here $D_{\RR^b}$ is the Gau\ss-Bonnet operator on $\RR^b$, $P$ is a self-adjoint  first order differential operator 
on $\Omega^*(F)$, and  
\begin{equation}\label{PA}
\begin{split}
(P+1/2)^2  = A_+, \\
(P-1/2)^2 = A_-.
\end{split}
\end{equation}
If $(M,g)$ is a simple edge space with an admissible edge metric, then $D$ is a perturbation of $(\del_s + s^{-1}P) + D_Y$ with 
higher order correction terms determined by the curvature of the Riemannian submersion $\phi: Y \to B$ and the second 
fundamental forms of the fibres $F$ in $Y$. 

Just as for the Hodge Laplacian, we can define the minimal and maximal extensions, $D_{\min}$ and $D_{\max}$ 
of $D$; their domains are $\calD(D_{\min}) = \calD_{\min}(D)$ and  $\calD(D_{\max}) = \calD_{\max}(D)$. We now 
characterize elements in these domains using arguments similar to those in Lemma~\ref{max}. 

\begin{lemma}\label{maxD}
Let $(M,g)$ be compact with simple edge singularity. Then, for any $\epsilon >0$
\begin{align*}
\calD_{\max}(D)\subset H^\infty_e\Omega^*(M) + x^{1-\epsilon} H^1_e\Omega^p(M) \subset H^1_e\Omega^*(M).
\end{align*}
Any $w \in \calD_{\max}(D)$ admits a weak asymptotic expansion 
\begin{align*}
&w \sim \sum_{j=1}^N x^{-\nu_j+1/2}  b_j(y,z)+\wt{w}, \quad \nu_j^2 \in \mbox{Spec}(A_+\oplus A_-) \cap (0,1), 
\ \nu_j > 0, \qquad \wt{w} \in x^{1-\epsilon} H^1_e\Omega^*.
\end{align*}
Note that $\nu_j=0$ does not appear in the expansion, since $x^{-1/2}\notin L^2(\R_+)$.

The coefficients $b_j(y, z)$ are smooth in $z$, but typically have negative Sobolev regularity in $y$. 

On the other hand, 
\begin{align*}
\calD_{\min}(D)=\{w\in \calD_{\max}(D)\mid b_j(y,z)= 0, j=1,...,N\}.
\end{align*}
\end{lemma}
\begin{proof}
Consider $w=(w_+,w_-)\in  \calD_{\max}(D)$ with respect to the decomposition of $D$ into even and odd components. 
The existence and nature of the terms in the weak asymptotic expansion for elements $w_{\pm}$ 
is derived exactly as for the Laplacian in Lemma~\ref{max}, see \cite[Theorem 7.3]{Maz-edge}. 

The exponents $\gamma$ in the expansion of $w_+$ arise as indicial roots of $(s\partial_s+P)$, 
in the interval $(-1/2,1/2)$. Consequently $\gamma \in -\mbox{Spec}\,(P)\cap (-1/2,1/2)$, 
with the coefficient $b_+(-\gamma)\in \calD'(Y, E(-\gamma))$ being a 
distribution in $Y$ with values in the finite dimensional $(-\gamma)$-eigenspace of $P$. 
Since on forms of even degree, $(P+1/2)^2=A_+$, we see that $\gamma=-\nu+1/2, \ \nu^2 \in \mbox{Spec}A_+ \cap (0,1)$. 
Hence the expansion of $w_+$ takes the form
\begin{align*}
w_+\sim \sum x^{\gamma}b_+(-\gamma) + \wt{w}_+, \
\wt{w}_+ \in x H^1_e\Omega^*.
\end{align*}

Similarly, the exponents $\gamma$ in the expansion of $w_-$ arise as indicial roots of $(-s\partial_s+P)$, 
in the interval $(-1/2,1/2)$. Consequently $\gamma \in \mbox{Spec}\,(P)\cap (-1/2,1/2)$, 
with the coefficient $b_-(\gamma)\in \calD'(Y, E(\gamma))$ being a distribution in $Y$ with values 
in the finite dimensional $\gamma$-eigenspace of $P$.
Since on the odd components $(P-1/2)^2=A_-$, we find again $\gamma=-\nu+1/2, \ \nu^2 \in \mbox{Spec}A_- \cap (0,1)$.
Hence the expansion of $w_-$ has the form
\begin{align*}
w_-\sim \sum x^{\gamma}b_-(\gamma) + \wt{w}_-, \
\wt{w}_- \in x H^1_e\Omega^*.
\end{align*}

It remains to establish the claim on $D_{\min}$. Note that if $\beta_{\pm}(\gamma) \in C^{\infty}(Y, E(\gamma))$ 
with $\gamma \in \mbox{Spec}\,(P)\cap (-1/2,1/2)$, then $(x^{-\gamma}\beta_+(\gamma), x^{\gamma}\beta_-(\gamma)\in \calD_{\max}(D)$.  
Now consider $w=(w_+,w_-)\in \calD_{\min}(D)\subset \calD_{\max}(D)$ with asymptotic expansion as above, and suppose that 
$v\in \calD_{\max}(D)$ is given by a sum of test sections 
$(x^{-\gamma}\beta_+(\gamma), x^{\gamma}\beta_-(\gamma))$, cut off smoothly at $x=1$ and extended to all of $M$. Then
\begin{align}
\langle Dw,v\rangle - \langle w,D v\rangle = \sum_{\gamma}
\langle b_+(\gamma), \beta_-(\gamma) \rangle + 
\langle b_-(\gamma), \beta_+(\gamma) \rangle,
\end{align}
where the summation is over $\gamma \in \mbox{Spec}\,(P)\cap (-1/2,1/2)$. This vanishes because $w$ and $v$ 
lie in adjoint domains, namely $\calD_{\min}$ and $\calD_{\max}$, respectively. Since the coefficients
$\beta_{\pm}(\gamma)$ are arbitrary, we deduce that each $b_{\pm}(\gamma) = 0$ as a distribution.  
Hence the vanishing of these
conditions is at least a necessary criterion to lie in the minimal domain. It is also sufficient since
any $v \in \calD(D_{\max})$ can locally be approximated in the graph norm by polyhomogeneous forms within the 
maximal domain, and if all coefficients for $w$ vanish, then $\langle Dw, v \rangle = \langle w, Dv \rangle$, and this is
sufficient to show that $w$ lies in the domain adjoint to $D_{\max}$, i.e.\ in $\calD(D_{\min})$.

The fact that we can approximate any $v \in \calD(D_{\max})$ locally in the graph norm by polyhomogeneous forms
follows by a standard mollification argument. First note that $v \in \calD(D_{\max})$ if and only if $\chi v \in 
\calD(D_{\max})$ for any smooth cutoff $\chi$ that is constant in $x$-direction near $x=0$; 
in other words, it suffices to work locally by using a partition of unity to reduce to functions compactly supported in product neighbourhoods $(x,y,z) \in (0,x_0) \times \calU \times F$, $\calU \subset \RR^b$. Now fix $\psi \in \calC^\infty_0(\RR^b)$ with $\int \psi =1$ and define
$\psi_\epsilon(y) = \psi_\epsilon(y):=\epsilon^{-b}\psi(y/\epsilon)$ so that $\int \psi_\e = 1$ for all $\e \to 0$.
By definition of the weak expansion, the convolution (in the $y$ variable only) 
$v_\epsilon = v \star \psi_\epsilon$ is polyhomogeneous in the strong sense. We also write $\Phi_\epsilon v = v_\epsilon$.
For differential forms we convolve the local coefficients of the form with $\psi_\epsilon$. We claim that $v_\e \in \calD(D_{\max})$ 
and $v_\e \to v$ in the graph norm. The fact that $v_\e \to v$ in $L^2$ is standard; furthermore,
$D v_\e = Dv \star \psi_\e + [D, \Phi_\epsilon] v$ lies in $L^2$, cf. \cite{BDV}, thus substantiating that $v_\e \in \calD(D_{\max})$,
and $D v_\e \to Dv$, which gives the graph convergence. 
\end{proof}

Finally, recall the characterization of the Friedrichs extension for the Hodge Laplacian
\begin{align*}
\Delta^{\mathrm{Fr}}=D_{\max}D_{\min}, \quad \calD (\Delta^{\mathrm{Fr}})=
\{w\in \calD_{\min}(D): Dw\in \calD_{\max}(D)\}.
\end{align*}
Thus, using the notation of Lemma \ref{max}, Lemma \ref{maxD} implies the following.
\begin{prop}\label{friedrichs-domain}
Let $(M,g)$ be compact with simple edge singularity and an adapted incomplete edge metric $g$, and 
$\Delta^{\mathrm{Fr}}$ the Friedrichs extension for the Hodge Laplacian. Then
$$ 
\calD(\Delta^{\mathrm{Fr}})=\{w\in \calD_{\max}(\Delta) \mid a^-_j(w;z)=0, j=1,...,N\}. 
$$
\end{prop}
This result is well-known especially in the conic case. Henceforth we omit the superscript $\mathrm{Fr}$
since we deal exclusively with the Friedrichs extension. Note that though we performed the discussion 
under rescalings $\Phi_p$ in each local coordinate neighborhood, index sets of rescaled solutions in 
$\calD_{\max}(\Delta)$ are invariant under coordinate changes, so the characterization of 
$\calD(\Delta^{\mathrm{Fr}})$ in Proposition \ref{friedrichs-domain} is globally well-defined.

\section{Construction of the heat kernel}\label{construct}
In this section we develop the `heat calculus' on a simple edge space $M$ and prove that the heat kernel of the 
Friedrichs extension of the Hodge Laplacian lies in this calculus. This construction follows a now standard path, 
and in fact is almost identical to the construction in \cite{Moo}  of the heat kernel for the Laplacian on spaces 
with isolated conic singularities; other related constructions appear in \cite{Mel-APS}, \cite{Albin}, and see also 
\cite{AM} for an expository account of several other heat kernel constructions using these same techniques of 
geometric microlocal analysis. 

The heat kernel is a priori a distribution on $[0,\infty) \times (\wt{M})^2$ which is $\calC^\infty$ on
$(0,\infty) \times (\mbox{int}\,M)^2$. The general idea is to study its
singular structure by lifting it to a `heat space' $M^2_h$, which is a  resolution of $\RR^+ \times (\wt{M})^2$ 
obtained by blowing up certain submanifolds of the boundary of this space; elements of the heat calculus are, by definition,
operators with Schwartz kernels which are polyhomogeneous on $M^2_h$. An iterative parametrix construction is 
used to construct a good approximation to the fundamental solution operator for $\del_t + \Delta_g$ in this class of
operators, and a standard regularity argument shows that the true heat kernel itself lies in this calculus too. 
Careful bookkeeping yields very precise information about the small-time asymptotics of this kernel near the 
singular strata, which leads eventually to the main result of this paper.

\subsection{The heat double space}
Fix an adapted coordinate chart $(x,y,z)$ on $M$ near the edge, or equivalently, on $\wt{M}$ near
its boundary. Taking two copies of this chart as well as the time coordinate $t$ yields a coordinate system
$(t,x,y,z,\tilde{x},\tilde{y},\tilde{z})$ on $\RR^+ \times (\wt{M})^2$, valid near the diagonal $(x,y,z)=(\tilde{x},\tilde{y},\tilde{z})$.
We define the resolution in two steps. The first is to blow up the fibre diagonal of $\partial \wt{M}$ at $t=0$, i.e.\ 
the submanifold
\[
\FD = \{ (0,0,y,z,0,\tilde{y},\tilde{z}): y = \tilde{y} \},
\]
using the parabolic homogeneity of the problem. The space which results from this is denoted $[ \RR^+; \FD, \{dt\} ]$; this
notation indicates both the submanifold that is to be blown up ($\FD$) and the direction $dt$ which is scaled parabolically.
We refer to \cite{EMM}; see also \cite{AM}, for a more careful discussion of parabolic blowups. This space is pictured in
the Fig. 2. 

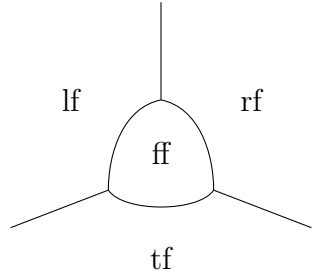
\begin{figure}[h]
\begin{center}
\begin{tikzpicture}
\draw (0,0.7) -- (0,2);
\draw (-0.7,-0.5) -- (-2,-1);
\draw (0.7,-0.5) -- (2,-1);

\draw (0,0.7) .. controls (-0.5,0.6) and (-0.7,0) .. (-0.7,-0.5);
\draw (0,0.7) .. controls (0.5,0.6) and (0.7,0) .. (0.7,-0.5);
\draw (-0.7,-0.5) .. controls (-0.5,-0.8) and (0.5,-0.8) .. (0.7,-0.5);

\node at (1.2,0.7) {\large{rf}};
\node at (-1.2,0.7) {\large{lf}};
\node at (0, -1.4) {\large{tf}};
\node at (0,0) {\large{ff}};
\end{tikzpicture}
\end{center}
\label{blowup-B}
\caption{The first blowup of the heat space.}
\end{figure}

It has four codimension one boundary faces, which we denote $\ff$, $\lf$, $\rf$ and $\tf$ (for `front face', `left face', 
`right face' and `temporal face'), respectively. One way of understanding this blowup is that it is the smallest manifold
with corners on which the lifts of smooth functions from $\RR^+ \times (\wt{M})^2$ and the `parabolic polar coordinates' 
$R = \sqrt{x^2 + \tilde{x}^2 + t}$, $\Theta = (t/R^2, x/R, \tilde{x}/R)$ are all $\calC^\infty$.  These polar
coordinates are difficult to use in computations, so we typically work with projective coordinate charts instead. 
Thus away from $\tf$ we use
\begin{equation}
\rho=\sqrt{t}, \  \xi=\frac{x}{\sqrt{t}}, \ \wt{\xi}=\frac{\wt{x}}{\sqrt{t}}, \ w=\frac{y-\widetilde{y}}{\sqrt{t}}, \ 
z, \ \widetilde{y}, \ \widetilde{z}.
\label{top-coord}
\end{equation}
Here $\rho, \xi$ and $\wt{\xi}$ are defining functions for the faces $\ff$, $\rf$ and $\lf$ away from $\tf$.  
On the other hand, away from $\lf$ we use the chart
\begin{equation}
\tau=\frac{t}{\widetilde{x}^2}, \ s=\frac{x}{\widetilde{x}}, \ \w=\frac{y-\widetilde{y}}{\widetilde{x}}, \ 
z, \ \widetilde{x}, \ \widetilde{y}, \ \widetilde{z},
\label{right-coord}
\end{equation}
so that $\tau, s, \widetilde{x}$ are defining functions of $\tf$, $\rf$ and $\ff$, respectively.  There are also analogous coordinates
valid away from $\rf$, obtained by interchanging the roles of $x$ and $\wt{x}$.  

Now, let $D_0 = \{(0,x,y,z,x,y,z)\}$ be the diagonal of $(\wt{M})^2$ at $t=0$ and $D_1$ its lift to this intermediate
heat space $[\RR^+ \times (\wt{M})^2; \FD, \{dt\}]$. In the coordinate system \eqref{right-coord}
above, 
\[
D_1 = \{\tau = 0, s=1, \omega = 0, z = \wt{z}\},
\]
respectively. The second and final step is to blow up $[\RR^+ \times (\wt{M})^2; \FD, \{dt\}]$ along $D_1$,
again parabolically in the $t$ direction.  This yields the final heat space:
\[
M^2_h = \left[ \big[ \RR^+\times (\wt{M})^2; \FD, \{dt\} \big]; D_1, \{dt\} \right]
\]
This has the same four boundary hypersurfaces as before, but also a new one, denoted $\td$ (for temporal diagonal), created in
the final blowup. This is illustrated in Fig. 3.
\begin{figure}[h]
\begin{center}
\begin{tikzpicture}
\draw  (0,0.7) -- (0,2);
\draw  (-0.7,-0.5) -- (-2,-1);
\draw  (0.7,-0.5) -- (2,-1);

\draw  (0,0.7) .. controls (-0.5,0.6) and (-0.7,0) .. (-0.7,-0.5);
\draw  (0,0.7) .. controls (0.5,0.6) and (0.7,0) .. (0.7,-0.5);
\draw  (-0.7,-0.5) .. controls (-0.5,-0.6) and (-0.4,-0.7) .. (-0.3,-0.7);
\draw  (0.7,-0.5) .. controls (0.5,-0.6) and (0.4,-0.7) .. (0.3,-0.7);

\draw  (-0.3,-0.7) .. controls (-0.3,-0.3) and (0.3,-0.3) .. (0.3,-0.7);
\draw  (-0.3,-1.4) .. controls (-0.3,-1) and (0.3,-1) .. (0.3,-1.4);

\draw  (0.3,-0.7) -- (0.3,-1.4);
\draw  (-0.3,-0.7) -- (-0.3,-1.4);

\node at (1.2,0.7) {\large{rf}};
\node at (-1.2,0.7) {\large{lf}};
\node at (1.1, -1.2) {\large{tf}};
\node at (-1.1, -1.2) {\large{tf}};
\node at (0, -1.7) {\large{td}};
\node at (0,0.1) {\large{ff}};
\end{tikzpicture}
\end{center}
\label{blowup-BD}
\caption{The blown up heat-space $M^2_h$.}
\end{figure}
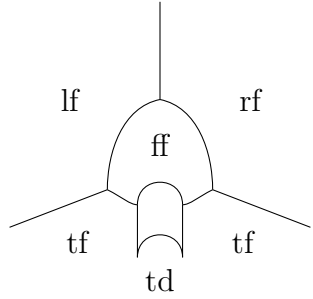

Using two pairs of interior coordinates $(x,y,z)$ and $(\wt{x},\wt{y},\wt{z})$ on $M$, then projective coordinates near $\td$ 
and away from $\ff$ are given by
\begin{equation}
T = \sqrt{t}, \ \Xi = ( (x-\wt{x})/T, (y-\wt{y})/T, (z-\wt{z})/T), \ \wt{x}, \wt{y}, \wt{z}.
\label{d-coord0}
\end{equation}
For a coordinate chart near $\td$ valid up to the corner $\td \cap \ff$ but away from $\tf$ and $\lf$, use the previous 
projective coordinate system to define
\begin{equation}
\eta=\sqrt{\tau}, \ \sigma =\frac{s-1}{\eta},\ \zeta =\frac{z-\widetilde{z}}{\eta}, \ \mu=\frac{\w}{\eta}, \ \widetilde{x}, \ 
\widetilde{y}, \ \widetilde{z}.
\label{d-coord}
\end{equation}
Thus $T$ and $\eta$ are the defining functions for $\td$ in these latter two coordinate systems, while $\wt{x}$ is a
defining function for $\ff$ in the first of these projective coordinate systems. Away from $\td$, it is sufficient to
use the projective coordinate systems from the intermediate space. 

There is a canonical blow-down map
\[
\beta: \wt{M}^2_h \longrightarrow \RR^+ \times \wt{M}^2
\]
which can be written explicitly in any of the local coordinate systems above. It is pertinent that the restriction of
$\beta$ to $\td$ is a fibration, with each fibre $\beta^{-1}( (0,\w,\w) )$ a closed `parabolic' hemisphere $S^m_+$ 
and base the lifted diagonal of $(\wt{M})^2$, which is diffeomorphic to a copy of $\wt{M}$. Similarly, the restriction 
of $\beta$ to $\ff$ is a fibration over the diagonal of $\del \wt{M} \times \del \wt{M}$ with fibre equal to a `parabolic' 
quarter-sphere $S^{m+1}_{++}$. 

\subsection{Heat calculus and parametrix construction}
We now define the heat calculus on $(M,g)$. This is a set of operators, acting on differential $p$-forms, 
characterized through the polyhomogeneity properties of 
the lifts of their Schwartz kernels to the double heat space $M^2_h$, as sections of ${}^{ie} \Lambda^p(M) \boxtimes {}^{ie} \Lambda^p(M)$. 
Recall that a kernel $K_A$ on $\RR^+ \times \wt{M}^2$ acts as an operator in one of two ways: first, it carries sections of 
${}^{ie} \Lambda^p(M)$ over $\wt{M}$ to sections of the same bundle over $\RR^+ \times \wt{M}$ by
\[
\phi \longmapsto  \int_{\wt{M}} K_A(t,z,\tilde{z}) \phi(\tilde{z})\, d\tilde{z} , 
\]
and second, it acts on sections of ${}^{ie} \Lambda^p(M)$ over $\RR^+ \times  \wt{M}$ by
\[
f \longmapsto \int_0^t \int_{\wt{M}} K_A(t-s,z,\tilde{z}) f(s,\tilde{z})\, d\tilde{z}. 
\]
We continue with the definition of the heat calculus under the rescaling transformation $\Phi$,
which leads to uniform orders of asymptotic expansions in each degree $p$.

\begin{defn}\label{heat-calculus}  
Let $\calE = (E_{\lf}, E_{\rf})$ be an index family for the two side faces of $M^2_h$. We define $\Psi^{\ell,p,\calE}_{\eh}(M)$ 
to be the space of all (rescaled) operators $A$ with Schwartz kernels $K_A$ which are pushforwards from polyhomogeneous functions 
$\wt{K}_A$ on $M^2_h$, with index family $\{(-b-3+\ell +j,0): j \in \mathbb N_0\}$ at $\ff$, $\{ (-m + p + j, 0) : j \in \mathbb N_0 \}$ 
at $\td$, $\emptyset$ at $\tf$ and $\calE$ for the two side faces of $M^2_h$. When $p = \infty$, $E_{\td} = \emptyset$. 
For simplicity we usually denote the lifted Schwartz kernel simply by $K_A$ again. 
\end{defn}
This is called a calculus because, at least under certain restrictions on the index sets at the side faces, one has that
\[
\Psi^{\ell,p,\calE}_{\eh}(M) \circ \Psi^{\ell',p',\calE'}_{\eh}(M) \subset \Psi^{\ell+\ell', p + p', \calE''}_{\eh}(M),
\]
where $\calE''$ is a new index family constructed explicitly from $\calE$  and $\calE'$. We do not need the most general form of 
this composition rule, but only the special case when $p = p' = \infty$, which is easier to prove. We defer the rather
technical proof to the Appendix. 

To simplify notation in this section, we define the heat calculus and describe the parametrix construction 
only for scalar operators.  Viewing the kernels as sections of ${}^{ie} \Lambda^p(M) \boxtimes {}^{ie} \Lambda^p(M)$,
then the heat kernel construction on differential forms is exactly the same and the leading orders of asymptotics 
of the kernels as polyhomogeneous distributions on $M^2_h$ are the same, since the homogeneity of the delta 
function acting on ${}^{ie}\Omega^p(M)$ does not depend on the degree $p$ under the form-bundle ${}^{ie}\Lambda^pM$ 
trivialization. We leave details of this to the reader, and shall reintroduce the bundle notation when it is needed more specifically.

The main point of this section is that it is possible to construct a parametrix for the solution operator of 
$\calL = \del_t + \Delta$ as an element in $\Psi^{2,0,\calE}_{\eh}(M)$ for some index family $\calE$; indeed, 
using this parametrix one can then show that the exact solution operator itself lies in this same space too. 
We are only interested here in the case where $\Delta = \Delta_p^{\mathrm{Fr}}$, but the same construction applies 
equally well for the heat equation associated to many other natural geometric operators associated to an incomplete 
edge metric on $M$. This parametrix construction follows the scheme in \cite{Moo}, as explicated further in \cite{AM}. 

The construction is a standard one in geometric microlocal analysis: it proceeds by positing that the heat
kernel lies in $\Psi_{\eh}^{2,0,\calE}(M)$ for some $\calE$, and then solving for the successive terms in the expansion
at certain faces of $M^2_h$ using the restrictions to the these faces of the lift of $t\calL$. This leads
to a parametrix for which the remainder vanishes to high, or even infinite, order at these faces. A better
parametrix, which vanishes to infinite order at other faces, is constructed using the composition formula,
and a short argument then yields that the true heat kernel is an element of this heat calculus and that
the parametrix captures the form of the asymptotic expansions at all boundaries.

It is helpful to work with the operator $t\calL$ rather than just $\calL$, simply because it lifts to an operator
which is smooth at the faces $\ff$ and $\td$ of $M^2_h$; this change is unimportant since $\calL H = 0$ implies $t\calL H = 0$, etc.
To commence, then, consider the lift of $t \calL$ from $\RR^+ \times M^2$, where $\Delta$ acts on the left copy of 
$M$, to $M^2_h$. Thus locally near $\ff$, $\td$ and $\tf$, $t\calL$ is a $b$-operator on $M^2_h$, i.e.\  a sum of products of 
elements of $\calV_b(M^2_h)$. (However, near $\rf$ it can be put into this form only by multiplying by $r^2$). 
This is easily checked using the local projective coordinate systems above. For example, near $\ff$ we compute that 
$t\del_x = \rho \del_\xi$, $t\del_x^2 = \del_\xi^2$, and so forth. 
Because of this structure, it makes sense to restrict $t\calL$ to one of these boundary faces.  We are particularly
interested in its restrictions to $\td$ and $\ff$, which we call the normal operators of $t\calL$ at these faces,
and denote by $N_{\td}(t\calL)$ and $N_{\ff}(t\calL)$, respectively. Using the projective coordinate system $(T,\Xi,\wt{\xi})$
near $\td$, we compute that
\begin{equation}
N_{\td}(t\calL) = \frac12 T \del_T + \Delta_\Xi - \frac12 \sum \Xi_j \del_{\Xi_j},
\label{ntdL}
\end{equation}
while using the projective coordinates $(\tau,s,w,z,\wt{x},\wt{y},\wt{z})$ near $\ff$ we see that
\begin{equation}
N_{\ff}(t\calL )= \tau\left(\del_\tau -\partial_s^2 + s^{-2}(A-1/4) + \Delta^{\R^b}_{\w}\right),
\label{nffL}
\end{equation}
where $A$ is the rescaled operator \eqref{a} acting in the variable $z$. Note that $N_{\ff}(t\calL)$ acts tangentially 
to the fibres $S^{n+1}_{++}$ of $\ff$, i.e.\ it involves no derivatives with respect to $(\wt{x},\wt{y},\wt{z})$.

We now search for an element $H \in \Psi^{2,0,\calE}_{\eh}(M)$ for which $t\calL \circ H$ vanishes (to as high order
as possible) at as many faces of $M^2_h$ as can be arranged. The index family $\calE$ will be determined 
in the course of the construction. 

The expansion of $H$ near $\td$ has a `universal' form.  Indeed, the expression for $N_{\td}(t\calL)$ above does
not depend on the edge geometry, but is just the same as for the corresponding operator on a closed
manifold without edges or other singularities. Using the projective coordinates $T = \sqrt{t}$ and $\Xi = (\xi - \wt{\xi})/T$,
write $H \sim \sum T^{-m+j} H_j$ near this face. The reason we have chosen the leading exponent to equal $-m$ is to 
be compatible with the initial condition 
\[
\left. H \right|_{t=0} = \delta(\xi - \wt{\xi}) = T^{-m} \delta(\Xi).
\]
Applying $N_{\td}(t\calL)$ to this formal expansion for $H$ shows that 
\[
\left(\Delta_\Xi - \frac12 \Xi \cdot \del_\Xi + \frac{j-m}{2}\right) H_j = Q_j,
\]
where $Q_j$ is an inhomogeneous term depending on $H_0, \ldots, H_{j-1}$. In particular, $Q_0 = 0$
and $H_0 = (4\pi)^{-m/2}T^{-m}\exp(-|\Xi|^2)$. One can verify inductively that each $H_j$ decays rapidly
in $\Xi$, and moreover, that this construction is uniform in the parameter $\wt{\xi} = (\wt{x},\wt{y},\wt{z})$,
even up to the front face where $\wt{x} = 0$. We refer to \cite{Mel-APS} and \cite{AM} for more leisurely treatments 
of this step.

We can also solve away at least the leading order coefficient of $H$ at $\ff$. As in the previous case,
expanding $H$ near this face and writing $N_{\ff}(H)$ for its leading coefficient there, we see that 
\[
N_{\ff}(t\calL \circ H) = N_{\ff}(t \calL) \circ N_{\ff}(H).
\]
To make this vanish, it is reasonable to choose $N_{\ff}(H)$ to equal the fundamental solution for
the heat operator $N_{\ff}(t\calL)$.  This heat operator is simply $\tau$ times the heat operator
for the Friedrichs extension of the Hodge Laplacian on the product $C(F) \times \RR^b$. In other words, we set
\[
N_{\ff}(H) = H^{C(F)}(\tau,s,z,\widetilde{s}=1,\wt{z}) H^{\RR^b}(\tau, \w,\widetilde{\w}=0).
\]
The first factor here is the Friedrichs heat kernel of the rescaled operator 
\begin{align}\label{cf}
\Delta^{C(F)}=-\frac{d^2}{ds^2} + \frac{1}{s^2}\left(A-\frac{1}{4}\right),
\end{align}
as in \cite{Ch1}, \cite{Moo} and \cite{Lesch}.  Note that $N_{\ff}(H)$ acts on the fibres of $\ff$. 
It has a singularity at $\tau = 0, s=1, \w = 0$ in $\ff$, which is precisely the intersection
$\ff \cap \td$. This matches the singularity of $N_{\td}(H)$ along $\td$ at this same intersection; this follows
from the fact that $N_{\td}(t\calL) H_0 = 0$ has a unique tempered solution up to scale, and we have chosen
the normalizing constant to make the value at $T=0$ equal to the delta function. We refer to \cite{AM} for
more on this point. 

We now find an element $H^{(1)} \in \Psi^{2,0,\calE}_{\eh}(M)$ with this prescribed asymptotic data at these two faces.
It was already explained that because of the homogeneity of the delta function, the leading term at $\td$ has
homogeneity $-m$. There is a similar phenomenon at $\ff$: we are integrating with respect to the volume form
of an incomplete edge metric, which has a factor of $x^f$, but there is compensating factor of $x^{-m} = x^{-b-f-1}$.
One must check through the density factors in projective or polar coordinate systems at this face, cf.\ \cite{Moo} or 
\cite{AM}, but the effect is that the leading exponent at $\ff$ is $(-b-1)$; this is precisely the 
homogeneity order of the delta function acting on ${}^{ie}\Omega^p (M)$ under the rescaling transformation $\Phi$. 
Thus we fix a smooth cutoff function 
$\chi(a)$ which equals $1$ for $x \leq a \ll 1$ and which vanishes outside a slightly larger interval, and can then regard 
$\chi(\rho_{\ff}) \rho_{\ff}^{-1-b}N_{\ff}(H)$ as defined on $M^2_h$. Choose some function $H'$ near $\td$ which has the 
given sequence of functions $H_j$ as the coefficients in its asymptotic expansion as $T \to 0$. By the compatibility discussed above, 
we see that these two requirements for $H^{(1)}$ fit together at $\ff \cap \td$.  

This discussion provides most of the proof of the
\begin{prop}
There exists an element $H^{(1)} \in \Psi^{2,0,\calE}_{\eh}(M)$, where $E_{\lf}= E_{\rf}$ is an index set determined by the set 
of  indicial roots of $x^2 \Delta^\Phi$, i.e.\ the numbers $\frac12 + \nu_j$, $\nu_j^2 \in \mbox{spec}\,(A)$,
such that $t\calL H^{(1)} = P^{(1)} \in \Psi^{3,\infty,\calE^{(1)}}_{\eh}(M)$, where  $\calE^{(1)} = (E_{\lf}, E_{\rf}-1)$, and with
\[
\lim_{t\to 0}H^{(1)}(t, x, y, z, \wt{x},\wt{y},\wt{z}) = \delta(x-\wt{x})\delta(y-\wt{y})\delta(z-\wt{z}).
\]
\end{prop}
\begin{proof}
The only point which is not immediately apparent from the construction of $H^{(1)}$ above is the characterization of
its expansions at $\lf$ and $\rf$.  These are inherited directly from the expansions of the Friedrichs heat kernel 
on the cone, and so we turn to a closer examination of this kernel. 

The Hodge Laplacian on the complete cone $C(F)$ can be described using separation of variables.  Following \cite[Prop. 2.3.9]{Lesch},
and using the eigendecomposition of the fibre operator $A$, $\Delta^\Phi_p$ separates into the ordinary differential operators
\[
l_{\nu}:=-\frac{d^2\,}{ds^2} + \frac{1}{s^2}\left(\nu^2-\frac{1}{4}\right),
\]
where $\nu^2 \in \mbox{spec}\,(A)$. There is an explicit formula for the corresponding Friedrichs heat kernel:
\begin{equation}
e^{-\tau l_{\nu}}(s,\widetilde{s})=\frac{1}{2\tau }(s\widetilde{s})^{1/2}I_{\nu}\left(\frac{s\widetilde{s}}{2\tau}\right)e^{-\frac{s^2+\widetilde{s}^2}{4\tau }} 
\label{bs0}
\end{equation}
where $I_\nu$ denotes modified Bessel function (the Bessel function with imaginary argument) of order $\nu$. 
Thus the heat kernel on $C(F)$ is given by
\begin{equation}\label{bessel-sum}
H^{C(F)}(\tau, s, z, \widetilde{s},\widetilde{z})= \sum_\nu \frac{(s \widetilde{s})^{\frac12}} {2\tau} I_{\nu}\left(\frac{s\widetilde{s}}{2\tau}\right)
e^{-\frac{s^2+\widetilde{s}^2}{4\tau}}\phi_{\nu}(z)\phi_{\nu}(\widetilde{z}),
\end{equation}
where $\phi_{\nu}$ is the eigenform associated to $\nu^2\in \mbox{spec}\,(A)$, and finally,
\begin{equation}
H^{C(F) \times \RR^b}(\tau,s,z,\w,\widetilde{s},\widetilde{z},\widetilde{\w}) = H^{C(F)}(\tau,s,z,\widetilde{s},\widetilde{z})H^{\RR^b}(\tau,\w,\widetilde{\w}),
\end{equation}
where
\[
H^{\RR^b}(\tau, \w, \widetilde{\w})=\frac{1}{(4\pi)^{b/2}} \frac{1}{\tau^b} e^{-|\w-\tilde{\w}|^2/4\tau}
\]
is the standard Euclidean heat kernel on $\RR^b$. Classical bounds for the Bessel functions show that 
this sum converges locally uniformly in $\calC^\infty$.    

Now, $\RR^+ \times C(F) \times \RR^b$ is naturally identified with the interior of each fibre of $\ff$, e.g.\ using 
the projective coordinates \eqref{right-coord}, and in terms of these, $s = x/\widetilde{x}$, 
$\w = (y-\widetilde{y})/\widetilde{x}$, so we evaluate $\tilde{s} = 1$, $\widetilde{\w} = 0$. The 
leading coefficient on $\ff$ is $H^{C(F)\times \RR^b}(\tau, s, z, \w, \widetilde{z}, 0)$. 
Using this identification and following \cite{Moo}, we see that $H^{C(F)\times \RR^b}$ extends to be 
polyhomogeneous on the front face of $M^2_h$. Furthermore, 
the index sets of $H^{(1)}$ at the left and right faces are exactly the same as those for $H^{C(F)}$. The index set at $\rf$, i.e.\
as $s \to 0$, can be read off directly from \eqref{bessel-sum}; the exponents are simply the indicial roots $(\nu+1/2),\nu \geq 0$.
By symmetry of the heat kernel, the index set at $\lf$ is exactly the same.  On the other hand, the index set for the
remainder term $P^{(1)}$ at $\rf$ is lowered by one because the operator $t\calL$ lowers terms by two orders but the leading
term is killed. 
\end{proof}

The next step in this construction involves choosing a slightly finer parametrix $H^{(2)}$ which leaves an error
which also vanishes to infinite order along the entire right face as well.  This is done by an iterative construction
similar to, but even simpler than, the procedure above at $\td$. 
\begin{prop}
There exists an element $H^{(2)} \in \Psi^{2,0,\calE}_{\eh}(M)$, where $\calE$ is the same index family as above,
such that $t\calL H^{(2)} = P^{(2)} \in \Psi^{3,\infty, E_{\lf}, \infty}_{\eh}(M)$ and $\lim_{t \to 0} H^{(2)} = \mathrm{Id}$. 
\end{prop}
\begin{proof}
This step proceeds exactly as in \cite{Moo}. The error term $P^{(1)}$ from the previous step has an expansion 
along $\rf$, and we use $t\calL$ to subtract off the successive terms of this expansion. Let $J$ denote
a kernel which has asymptotic sum  at this face with terms equal to the ones obtained in this procedure;
note that we can assume that $J$ vanishes to first order at $\ff$ and to infinite order at $\tf$. We
see that $H^{(2)} = H^{(1)} + J$ has all the desired properties.

In order to eliminate a term $s^{\gamma}a$ in the asymptotic expansion of the error term $P^{(1)}$ at 
$\rf$, it is only necessary to solve the indicial equation
\[
 (-\partial_s^2 + s^{-2}(A-1/4)) u = s^{\gamma}\left(\tau^{-1}a\right). 
\]
This is because all other terms in the expansion of $t\mathcal{L}$ at $\rf$ lower the exponent in $s$ 
by at most one, while the indicial part lowers exponent by two.   Note that $\tau, \w, \wt{x}, \wt{y}$ 
and $\wt{z}$ only enter this equation as parameters.  We have already discussed how to solve this
equation on the cone $C(F)$ using the Mellin transform. The solution is polyhomogeneous in all
variables, including the parameters. It has leading order $\gamma+2$ at $\rf$.

Iterating this argument, asymptotically summing these solutions, and adding the resulting 
kernel to $H^{(1)}$, we obtain a new parametrix $H^{(2)}\in \Psi_{\eh}^{2,0,\mathcal{E}}(M)$ 
with the same index set $\mathcal{E}=(E_{\lf}, E_{\rf})$ at the right and left boundary faces, 
and with a new error term $P^{(2)}\in \Psi_{\eh}^{3,\infty, E_{\lf},\infty}(M)$.
\end{proof}

We have now constructed a parametrix $H^{(2)}$ such that $t\calL H^{(2)} = P^{(2)}$ vanishes to infinite order
along $\rf$ and all along the bottom faces $\tf$ and $\td$, and which vanishes to higher order at $\ff$.
We can regard any of these kernels as acting on functions on $\RR^+ \times M$ in the usual way, by
\[
(Ku)(t,x,y,z) = \int_0^t \int_M K(t-s,x,y,z,\wt{x},\wt{y},\wt{z}) u(s, \wt{x}, \wt{y}, \wt{z})\, ds d\wt{x} d\wt{y} d\wt{z}.
\]
The identity operator corresponds to the kernel $K_{\mathrm{Id}} = \delta(x - \wt{x})\delta(y - \wt{y}) \delta(z - \wt{z})$. Viewed
in this way, 
\[
\calL H^{(2)} = \mathrm{Id} + \frac{1}{t} P^{(2)}. 
\]
Observe that $t^{-1} P^{(2)} \in \Psi^{2,\infty, E_{\lf}, \infty}_{\eh}$. 
The final stage in the parametrix construction is to consider the formal Neumann series 
\[
(\mathrm{Id} + P^{(2)})^{-1} = \mbox{Id} + \sum_{j=1}^\infty (- t^{-1} P^{(2)})^j := \mathrm{Id} + P^{(3)},
\]
Note here that composition formula for the heat edge calculus, which we discuss in the Appendix, yields 
$t^{-j} (P^{(2)})^j \in \Psi^{2j,\infty, E_{\lf},\infty}_{\eh}$.
In other words, successive composition of $t^{-1} P^{(2)}$ with itself produces an operator which vanishes to higher
and higher order on $\ff$. 

We then define
\[
H^{(3)} = H^{(2)}(\mathrm{Id} + P^{(3)}).
\]
In fact, a slightly finer analysis, see \cite{Mel-APS} and also \cite{BGV}, shows that this formal series
is actually convergent. We refer to these sources for the necessary estimates, both in the case of compact manifolds
and for manifolds with cylindrical ends, but it is easily seen that everything there transports immediately
to this setting (and many other ones as well). Indeed, this is a general feature of such Volterra series.

The final operator $H^{(3)}\in \Psi^{2,0,E_{\lf}, E_{\rf}}_{\eh}$ satisfies $t\calL H^{(3)} = \mathrm{Id}$ in the
sense of operators on $\RR^+ \times M$; equivalently, $t\calL H^{(3)} = 0$ as an operator from functions
on $M$ to functions on $\RR^+ \times M$, and $H^{(3)} = \mathrm{Id}$ at $t=0$. 

We have now produced a kernel $H = H^{(3)}$ which has the following three properties:
\begin{itemize}
\item[$\bullet$] $t\calL H = 0$;
\item[$\bullet$] $\left. H \right|_{t=0} = \delta(x-\tilde{x})\delta(y-\tilde{y}) \delta(z-\tilde{z})$;
\item[$\bullet$] The range of $H$ lies in the Friedrichs domain of $\Delta$ for all $t > 0$ fixed. Indeed, 
by construction, the range of $H$ lies in $\dom_{\max}(\Delta)$ and the index sets $(E_\lf, E_\rf)$ of $H$ correspond to 
the characterization of the Friedrichs domain in Proposition \ref{friedrichs-domain}.
\end{itemize}
These three properties characterize the Friedrichs heat kernel uniquely. Hence $H = H^M_{\mathrm{Fr}}$.

We have proved the
\begin{prop}
\label{hkprop1}
Let $(M,g)$ be a compact simple edge space with asymptotically admissible metric $g$. Then the heat
kernel $H_p$ for the Friedrichs extension of the (rescaled) Hodge Laplace operator on $p$-forms on $M$ is an element of
the heat space $\Psi^{2,0,\calE}_{\eh}(M; {}^{ie}\Lambda^p M)$, where $\calE=(E_\lf, E_\rf)$ with
\[
E_{\lf} = E_{\rf} = \{ (\nu + \frac12 + \ell, p): \nu \geq 0, \nu^2 \in \mathrm{spec}\,(A),\ \ell, p \in \mathbb N_0\}.
\]
\end{prop}
Note that we assert the absense of logarithmic terms in the $\ff$ expansion of the heat kernel.

\subsection{The even subcalculus}

The parametrix construction we have given already contains a significant 
amount of information about the asymptotic properties of the heat kernel at each of the boundaries of $M^2_h$. 
We now explain how, if one makes stronger assumptions about the asymptotic structure of the metric $g$ near the 
edge, one obtains stronger conclusions about the structure of these asymptotics.  This will be important  in the 
next section when we examine the heat trace and its asymptotics. 

To be more specific, we show that if the asymptotically admissible metric $g$ is even, 
then the Friedrichs heat kernel lies in a distinguished subalgebra of the heat calculus which
we call the even subcalculus. Elements of this even subcalculus are characterized by certain 
parity conditions on the coefficients in their asymptotic expansions, near both the faces $\td$ 
and $\ff$ of $M^2_h$.  As we exploit later, because of these parity conditions, many terms in the 
heat trace expansions of such elements vanish. Conditions of this type appear in many places 
in the literature, see \cite{Scott} for example, for a discussion of the closely related  
Kontsevich-Vishik canonical trace on a similar subclass of pseudodifferential operators on a closed manifold. 
Even more closely related to the discussion here, Melrose defines an even subcalculus for heat operators
on closed manifolds and on manifolds with asymptotically cylindrical ends \cite{Mel-APS}, and we refer also
to the even subcalculi in the papers by Albin \cite[Cor. 8.9]{Albin} and Albin-Rochon \cite{AR} in other geometric settings.

We impose separate `evenness' conditions on the asymptotic expansions of elements of $\Psi_{\eh}^*$ at the
two faces $\td$ and $\ff$. The condition at $\td$ is exactly the one introduced in \cite{Mel-APS} (and which
appears at least implicitly in many other places), while the one at $\ff$ is directly inspired by this but tailored
to the specific geometry.  Rather than considering these conditions separately, we shall define the subclass
of operators which satisfy evenness conditions at both faces. 

\begin{defn}\label{def-even}
Let $K_A$ be an element of the heat calculus $\Psi^{\ell,p,\calE}_{\eh}(M)$.  Then $K_A$ is an element of
the {\rm even subcalculus} if the following conditions hold:
\begin{enumerate}
\item[i)] At the face $\td$, and in terms of the projective coordinates $(T, \Xi, \widetilde{\xi})$ (which are valid in the interior 
of this face), suppose that 
\[
K_A \sim T^{-m}\sum_{j=0}^\infty \kappa^{\td}_j(\Xi, \widetilde{\xi}) T^{j};
\]
then we require that
\[
\kappa^{\td}_{j}(-\Xi,\widetilde{\xi}) = (-1)^j \kappa^{\td}_{j}(\Xi,\widetilde{\xi})
\]
for all $j \geq 0$.
\item[ii)] At the front face $\ff$, and in the projective coordinates $(\rho,\xi,\widetilde{\xi}, \w, \widetilde{y}, z,\widetilde{z})$ 
from \eqref{top-coord}, valid in the interior of that face, suppose that 
\[
K_A \sim \rho^{-b-3}\sum\limits_{j=0}^{\infty} \kappa^{\ff}_{j}(\xi, \widetilde{\xi}, \w, \widetilde{y}, z, \widetilde{z}) \rho^{j};
\]
then we require that 
\[
\kappa^{\ff}_j(\xi,\widetilde{\xi}, -\w, \widetilde{y}, z, \widetilde{z})=(-1)^j\kappa^{\ff}_{j}(\xi,  \widetilde{\xi}, 
\w, \widetilde{y}, z, \widetilde{z}).
\]
\end{enumerate}
We shall denote the set of all operators which satisfy these conditions by $\Psi^{\ell,p,\calE}_{\eh; \mathrm{evn}}(M)$.
\label{evencalcdef}
\end{defn}

Moreover, we may define the odd subcalculus $\Psi^{\ell,p,\calE}_{\eh; \mathrm{odd}}(M)$ by requiring 
the coefficients $\kappa^{\ff}_j$ to be odd in $\w$, if $j$ is even, and vice versa. 
Both the even and odd heat calculus are invariant under even changes of special coordinates, since we
specify front face behaviour of all kernel components in exactly the same way. 

There are a few different motivations for this definition. First, a straightforward calculation, see \cite[\S 7.1]{Mel-APS}
for the computation at $\td$, shows that if $K_A \in \calC^\infty(\wt{M}\times \wt{M} \times [0,\infty))$, then its lift to $M^2_h$ 
satisfies both conditions i) and ii). Second, it is also not hard to check that the initial parametrix $H^{(1)}$ in the construction 
above can be chosen to satisfy both of these conditions.  On the other hand, it is not immediately obvious that $\Psi^*_{\eh;\mathrm{evn}}$ 
is actually a subalgebra, i.e.\ closed under composition, which gives substance to the following: 
\begin{prop}
\[
\Psi^{k,\ell,\calE}_{\eh;\mathrm{evn}} \circ \Psi^{k',\ell',\calE'}_{\eh;\mathrm{evn}} 
\subset \Psi^{k+k',\ell+\ell',\calE''}_{\eh;\mathrm{evn}}.
\]
\label{evencomp}
\end{prop}
The rather technical proof is in the Appendix, following the proof of composition for the full heat calculus.

To establish that the exact heat kernel $H$ lies within the even subcalculus for an 
even asymptotically admissible metric, we need to prove invariance of $\Psi^*_{\eh;\mathrm{evn}}(M)$
under $\Delta$ and $t\calL$. This requires a careful treatment of kernels as operators acting on $p-$forms.
We may certainly assume that all kernels to be compactly supported in a special coordinate chart 
near the edge and work in these coordinates, so long as we show that all properties are invariant
under change to another equivalent special coordinate chart.  However, the evenness criterion is stable 
under even changes of coordinates, so it suffices to establish invariance of $\Psi^*_{\eh;\mathrm{evn}}(M)$ 
under $\Delta$ and $t\calL$ for one choice of special coordinates, which we fix now. 

We call a differential operator $P$ \emph{even} if its lift to $M^2_h$ satisfies 
\[
\beta^* P  \circ \Psi^*_{\eh; \mathrm{evn}}(M) \subset \Psi^*_{\eh;\mathrm{evn}}(M), \quad 
\beta^* P  \circ \Psi^*_{\eh; \mathrm{odd}}(M) \subset \Psi^*_{\eh;\mathrm{odd}}(M),
\]
while $P$ is said to be \emph{odd} if
\[
\beta^* P  \Psi^*_{\eh; \mathrm{evn}}(M) \subset \Psi^*_{\eh;\mathrm{odd}}(M), \quad 
\beta^* P  \Psi^*_{\eh; \mathrm{odd}}(M) \subset \Psi^*_{\eh;\mathrm{even}}(M).
\]
To verify evenness of $\Delta$ and $t\calL$, it clearly suffices to check that the de Rham differential 
$d$ is an even operator and the Hodge star $\star$ is odd.  

To verify these we must choose to work either with respect to the rescaled form bundles
or the unrescaled ones. In order to use the former, it is necessary to show that the
rescaling operator $\Phi$ preserves parity, and we explain the sense in which this is true now.
If $K$ is any Schwartz kernel acting on the unrescaled form bundles, then the same operator
acting on the scaled form bundles is given by $K^\Phi = \Phi^{-1} K \Phi$.
For any $w\in \Lambda^p M$ there is a local decomposition 
\begin{align}\label{b-2}
w = w_1 \wedge \, dx + w_2, \quad  w_1 \in \Lambda^{p-1} Y, w_2 \in \Lambda^p Y.
\end{align}
Hence, each $w\in \Lambda^p M$ can locally be viewed as an element of 
$\Lambda^{p-1} Y \oplus \Lambda^p Y$. Moreover, we may decompose $\Lambda^p Y$ locally as follows
$$
\Lambda^p Y = \bigoplus_{i+j=p} \Lambda^i B \otimes \Lambda^j F =: \bigoplus_{i+j=p} \Lambda^{i,j} Y.
$$
In terms of this, the action of an operator $K$ mapping $\Lambda^{i,k-1} Y \oplus \Lambda^{i,k} Y$ to 
$\Lambda^{j,\ell-1} Y \oplus \Lambda^{j,\ell} Y$, then $K$ and $K^\Phi$ have the matrix forms
\begin{equation}
K=\left( \begin{array}{ll}
K_{\ell-1, k-1} & K_{\ell-1,k} \\ K_{\ell,k-1} & K_{\ell,k}
\end{array} \right), \qquad  K^\Phi=\left(
\begin{array}{ll}
K^\Phi_{\ell-1, k-1} & K^\Phi_{\ell-1 , k} \\ K^\Phi_{\ell ,k-1} & K^\Phi_{\ell,k}
\end{array} \right). 
\label{a-2}
\end{equation}
One calculates explicitly that 
\begin{equation}\label{a-4}
\begin{aligned}
K_{\ell-1,k-1} & = x^{\ell-1-f/2}(\wx)^{k-1-f/2} K^\Phi_{\ell-1,k-1}, \\ K_{\ell,k} & = x^{\ell-f/2}(\wx)^{k-f/2} K^\Phi_{\ell,k},
\end{aligned}
\quad
\begin{aligned}
K_{\ell , k-1} & = x^{\ell - f/2}(\wx)^{k-f/2-1} K^\Phi_{\ell,k-1}, \\
K_{\ell-1,k} & = x^{\ell-f/2-1}(\wx)^{k-f/2} K^\Phi_{\ell-1,k}.
\end{aligned}
\end{equation} 
From these we obtain
\[
\begin{aligned}
K_{\ell-1,k-1} & \sim \rho^{-b-3+\ell+k-2-f} \sum_{j=0}^\infty \kappa_{j,\ell-1,k-1}^{\ff} \, \rho^j, \\ 
K_{\ell,k} & \sim \rho^{-b-3+\ell+k-f} \sum_{j=0}^\infty \kappa_{j,\ell,k}^{\ff} \, \rho^j, 
\end{aligned} \quad 
\begin{aligned}
K_{\ell,k-1} & \sim \rho^{-b-3+\ell+k-1-f} \sum_{j=0}^\infty \kappa_{j, \ell, k-1}^{\ff} \, \rho^j, \\
K_{\ell-1,k} & \sim \rho^{-b-3+\ell+k-1-f} \sum_{j=0}^\infty \kappa_{j,\ell-1,k}^{\ff} \, \rho^j.
\end{aligned}
\]
All coefficients $\kappa_{*,*,*}^{\ff}$ are functions of $(\xi, \wt{\xi}, \w, \wy, z, \wz)$. 
Hence if $K^\Phi \in \Psi^*_{\eh;\mathrm{evn}}(M,\Lambda^p M, \Lambda^r M)$, where evenness is relative
to some choice of special coordinates, then the terms in the expansions for the components of $K$ satisfy 
\begin{align}\label{a-5}
\kappa_{j,*,*}^{\ff} (\xi, \wt{\xi}, -\w, \wy, z, \wz) = (-1)^j \kappa_{j,*,*}^{\ff} (\xi, \wt{\xi}, \w, \wy, z, \wz)
\end{align}
for any choice of $*$. Similar formul\ae\ hold if $K^\Phi \in  \Psi^*_{\eh;\mathrm{odd}}(M)$. The importance of
this is that these are the same parity formul\ae\ as for the components of $K^\Phi$ despite the fact that the 
prefactors, i.e.\ the powers of $\rho$ are different between different components, while they are all
the same (namely $\rho^{-b-3}$) for the components of $K^\Phi$.  However, this also means that if we verify 
that an operator is even or odd with respect to the unrescaled bundles, in the sense that the equations
\eqref{a-5} hold, then the corresponding parity conditions hold for the conjugated operator acting
between the scaled form bundles. In other words, it suffices to study the 
structure of $d$ and $\star$ on the unrescaled bundles without worrying about the parity of $\Phi$. 

We first examine how the splitting \eqref{a-2} changes under an even change of coordinates. After a short
calculation, one sees that the components of a form in $\Lambda^{i,k-1} Y \oplus \Lambda^{i,k} Y$ 
are transformed to the new coordinates
by a $2$-by-$2$ coordinate transformation matrix; the on-diagonal
entries of this matrix are even functions of $x$ and the off-diagonal components are odd. Note finally
that conjugation of an even or odd operator by a matrix of this form preserves parity. 

Next, with respect to the local splitting \eqref{b-2}, the exterior derivative $d$ can be written as
\begin{align*}
d=\left(
\begin{array}{cc}
d_{Y,k-1} & (-1)^{k-1}\partial_x \\ 0 & d_{Y,k}
\end{array}\right).
\end{align*}
If $\kappa^{\ff}_{j,*,*}$ satisfies \eqref{a-5}, then one has $\partial_\w \kappa^{\ff}_{j,*,*}(\cdot, -\w, \cdot)=(-1)^{j+1} 
\partial_\w \kappa^{\ff}_{j,*,*}(\cdot, \w, \cdot)$, and from this it is straightforward to check that $d$ is an even operator.
Similarly, if $g$ is an even asymptotically admissible metric, then the Hodge star operator $*_k$ mapping 
$\Lambda^{i,k-1} Y \oplus \Lambda^{i,k} Y$ to $\Lambda^{b-i,f-k} Y \oplus \Lambda^{b-i,f+1-k} Y$ takes the form
\begin{align*}
*_k=\left(
\begin{array}{cc}
A_{f-k,k-1}(x) & A_{f-k,k}(x) \\ A_{f+1-k,k-1}(x) & A_{f+1-k,k}(x)
\end{array}\right),
\end{align*}
where, by explicit calculation, $A_{f-k,k-1}(x)$, $A_{f+1-k,k}(x)$ are families of endomorphisms which are odd in $x$, 
while $A_{f+1-k,k-1}(x)$, $A_{f-k,k}(x)$ are endomorphisms depending evenly in $x$. It follows that $*_k$ is an odd operator.

We deduce from these calculations that $\Delta$ and $t\calL$ are even operators and hence preserve the even 
and the odd subcalculi. Following through the parametrix construction of the previous section we
conclude  the following result.
\begin{prop}
The fundamental solution $H$ for the heat equation associated to the Friedrichs extension of the (rescaled) Hodge Laplacian 
on a simple edge space $(M,g)$ with even asymptotically admissible metric is an element of 
$\Psi^{2,0,\calE}_{\eh;\mathrm{evn}}(M, {}^{ie}\Lambda^p M)$ for every choice of special coordinates. 
\label{hkeven}
\end{prop}

\section{Heat trace asymptotics}
Let us now turn to the form of the heat trace expansion, with particular attention to the implications of Proposition~\ref{hkeven}
on that expansion. 

First recall the diagonal $D_0 \subset \RR^+ \times \wt{M}^2$, and more importantly, its lift $D_h$ to $M^2_h$,
as pictured in Fig. 4.
\begin{figure}[h]
\begin{center}
\begin{tikzpicture}
\draw  (0,0) -- (0,2);

\draw  (0,0) .. controls (0.5,-0.5) and (1,-1) .. (1,-2);
\draw  (0,0) .. controls (0.5,-0.1) and (1.2,-0.5) .. (1.6,-1.5);
\draw  (0,0) .. controls (0.5,0) and (1.5,0) .. (2,-1.3);

\draw  (1,-2) -- (3,-3);
\draw  (2,-1.3) -- (4,-1.3);

\draw  (1.7,-1.5) .. controls (1.5,-1.5) and (1.4,-1.6) .. (1.4,-1.8);
\draw  (1.7,-1.5) -- (3.6,-1.9);
\draw  (1.4,-1.8) -- (3.4,-2.3);

\draw  (3.6,-1.9) .. controls (3.5,-1.9) and (3.4,-2.2) .. (3.4,-2.3);
\draw  (3.6,-1.9) .. controls (3.7,-1.9) and (3.8,-2) .. (3.8,-2.1);

\draw  (3.6,-1.9) -- (3.6,0.1);
\draw  (3.6,0.1) -- (0,0.7);

\draw  (1,-2) -- (1.4,-1.8);
\draw  (2,-1.3) -- (1.8,-1.5);

\node at (4.5,-2) {\large{td}};
\node at (-0.5,1) {\large{lf}};
\node at (0.1,-1) {\large{ff}};
\node at (4,1) {$D_h \subset M^2_h$};

\end{tikzpicture}
\end{center}
\label{diagonal-picture}
\caption{The diagonal hypersurface $D_h$ in $M^2_h$.}
\end{figure}
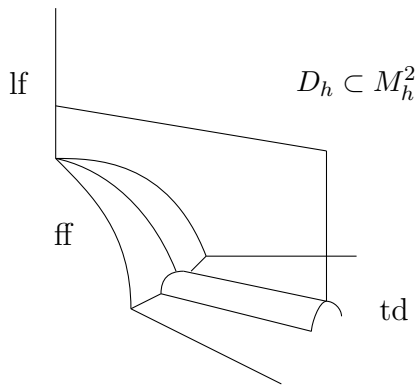

\noindent If $\iota_D: D_h \hookrightarrow M^2_h$ is the natural embedding and $H_k$ the heat
kernel for the Friedrichs extension of the Hodge Laplacian (which we consider as living on $M^2_h$
already, for convenience), then by Proposition \ref{hkprop1}, 
$\iota_D^* H_k$ is polyhomogeneous on $D_h$, and hence so is $\tr \iota_D^* H_k$, its pointwise trace. 
Our next goal is to obtain information about the terms which can appear in the expansion of this function. 
Let $D_{\td} = D_h \cap \td$, $D_{\ff} = D_h\cap \ff$ and $D_c = D_h \cap \rf\cap \lf$ be the
boundary faces of $D_h$. 

\begin{prop}\label{traceexpansion}
Let $g$ be asymptotically admissible. Then $\tr \iota_D^* H_k \in \calA_{\phg}^{\calG}(D_h)$,
where the index family is given by
\[
G_{\td} = -m + 2\mathbb N_0,\ G_{\ff} = -1-b + \mathbb N_0,\ G_c = E_{\lf} + E_{\rf}
\]
Moreover, if $g$ is even, then $G_{\ff} = -1-b + 2\mathbb N_0$. 
\end{prop}

This Proposition is an immediate consequence of the polyhomogeneous structure of $H_k$ itself. 
The fact that the restriction of a polyhomogeneous distribution to the submanifold $D_h$ 
is polyhomogeneous is straightforward and can be checked in local coordinates, but is also
a special case of the `pullback theorem' in \cite{Mel2}. 
Hence the form of the expansion in the general case follows directly from Proposition~\ref{hkprop1}, except
for the statement that the expansion at $\td$ has no odd terms.

The vanishing of the odd terms in the expansion at $\td$ in general and at $\ff$ when $g$ is even 
follows from Proposition~\ref{hkeven}. Indeed, the coefficients in the expansions of $\tr \iota_D^* H_k$
at each face of $D_h$ are the pointwise traces of the coefficients in the asymptotic expansions of $H_k$ 
at the corresponding faces of $M^2_h$. However, functions on $\td$ which are odd with
respect to the reflection $\Xi \mapsto -\Xi$, and similarly functions on $\ff$ which are odd with 
respect to $w \mapsto -w$, must vanish on $D_h$, and thus have vanishing trace. 
It is a general fact, explained in \cite{Mel-APS}, that in the interior, $H_k$ always satisfies the
evenness condition at $\td$, while this evenness at $\ff$ only holds if $g$ is even.

Our main interest, however, is not in the pointwise heat trace but in its integral, $\Tr H_k$.  
Note that $D_h$ can also be obtained as the parabolic blowup of $\RR^+ \times \wt{M}$
at $\{t = x = \wt{x} = 0,\ y = \wt{y}\}$ without any reference to the double space. Let
\[
\pi_c: D_h \to \RR^+
\]
be the composition of the blowdown map $D_h \to \RR^+ \times \widetilde{M}$ and the projection 
$\RR^+ \times \widetilde{M} \to \RR^+$. Then $\Tr H_k$ is the pushforward of $\tr \iota_D^* H_k$ by $\pi_c$:
\begin{equation}
\Tr H_k (t) = (\pi_c)_* \left( \tr \iota_D^* H_k \right) = \int_{\widetilde{M} \times \{t\}} \tr \iota_D^* H_k(t, x,y,z)\, x^{-f} dV_g.
\label{bigtrace}
\end{equation}
(The extra factor of $x^{-f}$ in the integrand is due to the rescaling $\Phi_k$.)
The formula with the integral on the right is valid when $t > 0$ since $D_h$ is the same as the diagonal $D_0$
of $\RR^+ \times \wt{M}^2$ there, but that this expression must be interpreted carefully near $t=0$
because of the front face of $D_h$. 

There is a good formalism for understanding pushforwards of polyhomogeneous distributions, known
as Melrose's Pushforward Theorem \cite{Mel2}. It states, roughly, that if $f: X \to Y$ is a special
type of map, known as a $b$-fibration, between two manifolds with corners, and if $u \in \calA_{\phg}^*(X)$,
then so long as the integrals along the fibres (or generalized fibres) make sense, we have
that $f_* (u) \in \calA_{\phg}^*(Y)$. Moreover, the index set of the image $f_*(u)$ is determined by
the index set of $u$, and similarly the coefficients in the expansions of the image can be determined
in terms of (in general, regularizations of similar sorts of pushforwards) coefficients of $u$. 
For the reader's convenience, we review this theorem here in this specific geometric setting.

We first note that the map $\pi_c$ is indeed a $b$-fibration. We refer to \cite{Mel2} for the precise
definition of this condition (see also the Appendix below), but do not repeat it here since it does not
play a significant role. Using the two sets of projective coordinates on $D_h$, arising 
from the coordinates \eqref{right-coord} and \eqref{d-coord} on $M^2_h$, we find
\begin{align}\label{pi}
\pi_c^*(t)=\rho_{\ff}^2\rho_{\td}^2.
\end{align}
It turns out that the formula for index sets of pushforwards is simpler 
if we write everything in terms of $b$-densities. Thus suppose that $\mu_0$ is a density
on $D_h$ which is smooth up to all boundary faces and everywhere nonvanishing. A smooth $b$-density 
$\mu_b$ is, by definition, any density of the form $(\rho_{\td} \rho_{\ff} \rho_c)^{-1}\mu_0$. There is
a pleasant naturality of $b$-densities with respect to blowups: a straightforward computation shows that
\[
\mu_b = A \beta_D^* ( (xt)^{-1} dt dxdydz)
\]
for some smooth nonvanishing function $A$. Note too that we can use $\mu_b = \beta_D^* ( t^{-1} x^{-f-1} dt dV_g)$
for this $b$-density. 

Consider a polyhomogeneous $b$-density $K$ on $D_h$, so $K = \tilde{K} \mu_b$ where $\tilde{K} \in \calA_{\phg}^*(D_h)$.
Now suppose that $\tilde{K}$ has index family $\calG = (G_\td, G_\ff, G_c)$. In order that the integrals
over the fibres of $\pi_c$ converge, one needs that all terms in the expansion of $\tilde{K}$ at $D_c$
have strictly positive exponent, i.e.\ $\mbox{Re}\, z > 0$ for all $(z,p) \in G_c$. Provided this integrability
condition holds, then the pushforward theorem asserts that 
\[
(\pi_c)_* K = (\pi_c)_* (\tilde{K} \mu_b) = T(t) \frac{dt}{t},\quad T \in \calA_{\phg}^F(\RR^+), 
\]
where the pushforward index set $F = \pi_c^{\flat}(\calG)$ is defined as follows. 
If $G_{\td} = \{ (z_j,p_j)\}$ and $G_\ff = \{w_\ell, q_\ell)\}$, then $G'_{\td} = \{ (z_j/2,p_j)\}$, $G'_\ff = \{(w_\ell/2, q_\ell)\}$
and $F=G'_\td \overline{\cup} G'_{\ff}$ is the `extended union' of $G'_{\td}$ and $G'_\ff$ 
\begin{align}
\label{extended}
G'_\td \overline{\cup} G'_{\ff} = G'_\td \cup G'_\ff \cup \{((z, p + q + 1): \exists \, (z,p) \in G'_\td,\ 
\mbox{and}\  (z,q) \in G'_{\ff} \}.
\end{align}
We refer to \cite[Eqn. (42)]{Mel2} for a proof. 
In our particular example we set $\tilde{K}=2\rho_c\rho_{\ff} \left( \tr \iota_D^* H_k \right)$, so that $T(t)=\Tr H_k (t)$;
and obtain the heat trace asymptotics as above.

We summarize all of this in the following central result.
\begin{thm}\label{heat-trace-corr}
Let $(M^m,g)$ be a compact simple edge space, where $g$ is an asymptotically admissible metric.
Let $H_k$ denote the heat kernel for the Friedrichs extension for the Hodge Laplacian on $k$-forms. 
Then, as $t \to 0$,
\[
\Tr H_k(t) \sim \sum_{\ell=0}^{\infty}A_\ell t^{\ell-\frac{m}{2}}+\sum_{\ell=0}^{\infty}C_\ell t^{\frac{\ell-b}{2}}+
\sum_{\ell \in \mathfrak{I}}G_\ell t^{\frac{\ell-b}{2}}\log t. 
\]
The index set $\mathfrak{I}$ here is 
\[
\mathfrak{I}=\{\ell\in \N_0 \mid \ell + m -b \ \mbox{even}\}.
\]

If the metric $g$ is even, then for every choice of special coordinates
\[
\Tr H_k(t) \sim \sum_{\ell=0}^{\infty}A_\ell t^{\ell-\frac{m}{2}}+
\sum_{\ell=0}^{\infty}C_\ell t^{\ell-\frac{b}{2}}+\sum_{\ell\in \mathfrak{I}'}G_\ell t^{\ell-\frac{b}{2}}\log t, 
\]
where $\mathfrak{I}' = \emptyset$ if $(m-b)$ is odd and $\mathfrak{I}'=\N_0$ if $(m-b)$ is even.
\end{thm}

We can say something about the coefficients which appear here. The coefficients
$b_{2\ell}$ of $\rho_{\td}^{-m+2\ell}$ in the expansion of $\tr \iota_D^* H_k$ (i.e.\ before integrating) are the standard local geometric 
quantities in the interior of $M$; on the other hand, the coefficients $a_{\ell}$ of $\rho_{\ff}^{-b-1+\ell}$ typically 
involve not only similar local quantities on the edge $B$ and on all of $M$, but also global quantities over $F$.
The coefficients $A_\ell$, $C_\ell$ and $G_\ell$ which appear above are integrals of these various coefficients which
involve both local and global quantities. We do not comment on 
this further since it plays no role here, but intend to return to a closer examination of these coefficients in 
another paper. We refer to \cite{Seeley} for an example where global spectral invariants on $F$ appear. 

We conclude this section by stating the implication of this main theorem on the meromorphic
structure of the zeta functions.  Let $P_k$ denote the projector onto the nullspace of the Friedrichs
extension of the Hodge Laplacian $\Delta_k$, and define
\[
\zeta_k(s):=\frac{1}{\Gamma(s)}\int_0^{\infty}t^{s-1}\Tr (H_k-P_k) \, dt, \quad \mbox{Re}\,(s) \gg 0.
\]
By the usual arguments, Theorem~\ref{heat-trace-corr} implies the following
\begin{prop}\label{zeta}
Each $\zeta_k(s)$ extends meromorphically to the entire complex plane. Near $s=0$, 
\[
\Gamma(s)\zeta_k (s) =
\sum_{j=0}^{\infty}\frac{A_j}{s+j-\frac{m}{2}}+\sum_{j=0}^{\infty}\frac{C_j}{s+\frac{j-b}{2}} - 
\sum_{j\in \mathfrak{I}}\frac{G_j}{\left(s+\frac{j-b}{2}\right)^2}, 
\]
where $\mathfrak{I}=\{j\in \N_0 \mid j + m -b \ \textup{even}\}$. If $g$ is also even, then
for every choice of special coordinates
\[
\Gamma(s)\zeta_k(s) = \sum_{j=0}^{\infty}\frac{A_j}{s+j-\frac{m}{2}}+
\sum_{j=0}^{\infty}\frac{C_j}{s+j-\frac{b}{2}} - \sum_{j\in \mathfrak{I}'}\frac{G_j}{\left(s+j-\frac{b}{2}\right)^2},
\]
where $\mathfrak{I}'=\emptyset$ if $(m-b)$ is odd and $\mathfrak{I}'=\N_0$ if $(m-b)$ is even.

In particular, if $m$ is odd, then $\zeta_k(s)$ is regular at $s=0$ 
for every $k$.  If $g$ is even and both $m$ and $b$ are odd, then 
for every choice of special coordinates we even have $\zeta_k(0) = 0$. 
\end{prop}

\section{Analytic torsion and the metric anomaly}\label{torsion}
We come, at least, to an examination of the analytic torsion of a compact simple edge space $(M,g)$ with an
asymptotically admissible edge metric. Recall that in terms of the zeta functions defined at the end of the last section, we define
the analytic torsion zeta function
\[
\zeta_{\AT}(s) := \frac12 \sum_{k=0}^m (-1)^k k \, \zeta(s,\Delta_k),
\]
and then set
\[
\log T(M,g) = \left. \frac{d\,}{ds}\right|_{s=0} \zeta_{\AT}(s).
\]

In certain cases, the analytic torsion zeta function $\zeta_{\AT}$ may be regular at $s=0$ even though
some or all of the individual zeta functions are not. For example, if $m$ is even, then there is nothing
in Proposition \ref{zeta} which precludes $\zeta (s, \Delta_k)$ from having a pole at $s=0$. In the special 
case of isolated conical singularities ($b=0$), it was shown in \cite{Dar} that despite this, there is a 
cancellation of residues and $\zeta_{\AT}(s)$ is regular at $s=0$. In fact, Dar's arguments also apply to 
more general situations where the conic cross-section $F$ may have boundary or even be singular. 
A consequence of our work here, however, is that if $F$ is smooth and $m$ is odd, then
each $\zeta_k(s)$ is regular at $s=0$, hence the analytic torsion is well-defined, without need for Dar's
cancellation argument. 

Recall from the introduction that we consider a slightly different quantity than just the analytic torsion.
The determinant line bundle over $M$,
\[
\det \mathcal{H}^*(M) :=\bigotimes_{k=0}^{m}\left(\, \bigwedge\nolimits^{m}\ker \Delta_{k}\right)^{(-1)^{k+1}},
\]
has a natural Hermitian inner product and norm $\|\cdot \|_{\, \det \mathcal{H}^*(M)}$ induced from the restriction
to $\ker \Delta_k$ of the $L^2$ norm on $L^2 \Omega^k$. The more invariant object is the Quillen metric
\[
\|\cdot \|^{RS}_{(M,g)}:= T(M,g)\|\cdot \|_{\, \det \mathcal{H}^*(M)}.
\]

We now study the variation of the analytic torsion on $M$ as we vary $g$. Consider a family $g_\mu$ of asymptotically 
admissible edge metrics on $M$,  $\mu \in (\mu_0 - \epsilon,\mu_0 + \epsilon)$. Denote by $\Delta_k(\mu)$ the Hodge 
Laplacian on $k$-forms associated to $g_\mu$.  In order to make these operators act on the same Hilbert space,
we conjugate by the natural local isometry
\[
T_{\mu} = \sqrt{\star_{\mu_0}^{-1} \star_{\mu}} \, \sqrt{\frac{dV_{g_{\mu}}}{dV_{g_{\mu_0}}}}: L^2\Omega^k(M,g_{\mu})\longrightarrow L^2\Omega^k(M,g_{\mu_0}).
\]
Thus consider the operators
$$
\wt{\Delta}_k(\mu):=T_{\mu} \circ \Delta_k(\mu) \circ T_{\mu}^{-1}
$$ 
as an operator on the fixed Hilbert space  $L^2\Omega^k(M,g_{\mu_0})$. 
The semigroup property of the heat kernel gives the identity 
\begin{equation}\label{difference}
\begin{split}
\frac{e^{-t\wt{\Delta}_k(\mu)}-e^{-t\wt{\Delta}_k(\mu_0)}}{\mu-\mu_0}&=
\int_0^t\frac{\partial\, }{\partial s}\left(\frac{e^{-(t-s)\wt{\Delta}_k(\mu_0)}e^{-s\wt{\Delta}_k(\mu)}}{\mu-\mu_0}\right)\, ds \\
&=\int_0^t e^{-(t-s)\wt{\Delta}_k(\mu_0)}\left(\frac{\wt{\Delta}_k(\mu_0)-\wt{\Delta}_k(\mu)}{\mu-\mu_0}\right)e^{-s\tilde{\Delta}_k(\mu)}\,ds.
\end{split}
\end{equation}
The identity \eqref{difference} is a relation between integral kernels, where $(\wt{\Delta}_k(\mu_0)-\wt{\Delta}_k(\mu))$
denotes a differential expression applied to the integral kernel of $\exp (-s\tilde{\Delta}_k(\mu))$.
Taking the limit $\mu \to \mu_0$ yields
\begin{align*}
\left. \frac{\partial\,}{\partial \mu}\Tr \left(e^{-t\wt{\Delta}_k(\mu)}\right)\right|_{\mu=\mu_0}=
&-\int_0^t\Tr \left(e^{-(t-s)\wt{\Delta}_k(\mu_0)}\left(\dot{\wt{\Delta}}_k(\mu_0)\right) e^{-s\wt{\Delta}_k(\mu_0)}\right)\\
=& -t \, \Tr\left(\left(\dot{\wt{\Delta}}_k(\mu_0)\right)  e^{-t\wt{\Delta}_k(\mu_0)}\right),
\end{align*}
where the upper dot denotes the derivative in $\mu$ of the differential expression $\wt{\Delta}_k$. 
Evaluating $\dot{\wt{\Delta}}_k(\mu)$ explicitly in terms of 
$T_{\mu}$, we find
\begin{multline*}
\left. \frac{\partial\, }{\partial \mu}  \Tr \left(e^{-t\wt{\Delta}_k(\mu)}\right) \right|_{\mu=\mu_0}= 
 -t \cdot \Tr \left(\dot{T}_{\mu_0}\circ \Delta_k(\mu_0)\circ e^{-t\Delta_k(\mu_0)}\circ T^{-1}_{\mu_0}\right) \\
-t \cdot \Tr \left(\Delta_k(\mu_0)\circ \dot{T}^{-1}_{\mu_0}T_{\mu_0}\circ e^{-t\Delta_k(\mu_0)}\right) 
 -t \cdot \Tr \left(\dot{\Delta}_k(\mu_0)\circ e^{-t\Delta_k(\mu_0)}\right).
\end{multline*}
Using that bounded operators can be commuted under the trace, the second term on the right here becomes
\begin{multline*}
\Tr \left(\Delta_k(\mu_0)\circ \dot{T}^{-1}T_{\mu_0}\circ e^{-t\Delta_k(\mu_0)}\right) 
= \Tr \left(e^{-t/2\Delta_k(\mu_0)}\Delta_k(\mu_0)\circ \dot{T}^{-1}_{\mu_0}T_{\mu_0}\circ e^{-t/2\Delta_k(\mu_0)}\right) \\
= \,  \Tr \left(\dot{T}^{-1}_{\mu_0}T_{\mu_0}\circ \Delta_k(\mu_0) \circ e^{-t \Delta_k(\mu_0)}\right) 
= - \Tr \left(T^{-1}_{\mu_0}\dot{T}_{\mu_0}\circ \Delta_k(\mu_0) \circ e^{-t\Delta_k(\mu_0)}\right). 
\end{multline*}
Consequently the first and second terms cancel and we are left with
\[
\left. \frac{\partial}{\partial \mu} \Tr \left(e^{-t\Delta_k(\mu)}\right)\right|_{\mu=\mu_0}=
 -t \, \Tr \left(\dot{\Delta}_k(\mu_0)\circ e^{-t\Delta_k(\mu_0)}\right).
\]
If $P_k(\mu)$ is the orthogonal projection onto the kernel of $\Delta_k(\mu)$, then repeating the arguments in
\cite[pp. 152-153]{RS}, we have
\begin{equation}\label{RS-identity}
\begin{split}
&\frac{\partial }{\partial \mu}\sum_{k=0}^m(-1)^k\cdot k \cdot \Tr \left(e^{-t\Delta_k(\mu)}-P_k(\mu)\right) \\
& \quad =t\, \frac{\partial }{\partial t}\sum_{k=0}^m(-1)^k\Tr \left(\A_{\mu}^k \left(e^{-t\Delta_k(\mu)}-P_k(\mu)\right)
\right)\, dt.
\end{split}
\end{equation}
Here $*_{\mu}$ is the Hodge star for $g_\mu$ and $\A^k_{\mu}:=*_{\mu}^{-1}\dot{*}_{\mu},$ on forms of degree $k$. Put
\begin{align*}
f(\mu,s):= \frac{1}{2}\sum_{k=0}^m (-1)^k\cdot k \cdot \frac{1}{\Gamma(s)}\int_0^{\infty}t^{s-1}\Tr 
\left(e^{-t\Delta_k(\mu)}-P_k(\mu)\right)\, dt,
\end{align*} 
so that, by definition, $\log T(M,g_{\mu})=\left.\frac{\partial}{\partial s}\right|_{s=0}f(\mu, s)$.

Now differentiate $f(\mu,s)$ under the integral. Inserting \eqref{RS-identity}, recalling that 
the heat trace decays exponentially and assuming $\mbox{Re}\,s \gg 0$ so that there are no
boundary terms in the integration by parts, we have
\begin{equation}\begin{split}\label{s-factor}
\frac{\partial }{\partial\mu}f(\mu, s)&=\frac{1}{2}\sum_{k=0}^{m}(-1)^{k}\frac{1}{\Gamma(s)}
\int_0^{\infty}t^{s}\frac{d}{dt}\Tr \left(\A_{\mu}^k \left(e^{-t\Delta_k(\mu)}-P_k(\mu)\right)\right)\, dt\\
&=\frac{1}{2}\, s\sum_{k=0}^{m}(-1)^{k+1}\frac{1}{\Gamma(s)}\int_0^{\infty}t^{s-1}\Tr 
\left(\A_{\mu}^k \left(e^{-t\Delta_k(\mu)}-P_k(\mu)\right)\right)\, dt.
\end{split}\end{equation}

At this point we introduce the separate hypotheses on the family of metrics $g_\mu$.  
First suppose that all the metrics are strongly asymptotically admissible, so that $g_\mu - g_{\mu_0} = 
\calO(x^{b+1})$. Then, there exists a family $\{e^j_\mu(p)\}$ of local orthonormal frames of $T_pM$, 
with $e^j_\mu(x, \cdot, \cdot) - e^j_{\mu_0}(x, \cdot, \cdot)= \calO(x^{b+1})$ for $p=(x,y,z)$. 
Consequently $\A_{\mu}^k(x)=\calO(x^{1+b})$ as well.  

Using projective coordinates, and letting 
$\pi_L:\R^+\times \widetilde{M}^2 \to \widetilde{M}$ be the projection onto the left copy of $\widetilde{M}$, 
we see that the lift $(\pi_L \circ \beta)^*\A^k_{\mu}$ is smooth on $M^2_h$ with 
\begin{equation}
\begin{split}
&\iota_D^* \tr (\A_{\mu}^k e^{-t\Delta_k(\mu)}) \sim \rho_{\ff}^0 \left(\sum\limits_{j=0}^{\infty}a_{j}\rho_{\ff}^{j}\right) \\ 
&\iota_D^*\tr (\A_{\mu}^k e^{-t\Delta_k(\mu)}) \sim \rho_{\td}^{-m}\left(\sum\limits_{j=0}^{\infty}b_{2j}\rho_{\td}^{2j}\right).
\end{split}
\end{equation}
Using the pushforward theorem once again, we deduce that
\[
\Tr \left(\A^k_{\mu}e^{-t\Delta_k(\mu)}\right) \sim \sum_{j=0}^{\infty}A_jt^{j-\frac{m}{2}}+\sum_{j=1}^{\infty}C_jt^{\frac{j}{2}}+
\sum_{j=1}^{\infty}G_jt^{\frac{j}{2}}\log t.
\]

On the other hand, suppose that each $g_\mu$ is even in $x$.  Then $\A_{\mu}^k(x)$ is also even, and hence
lies in the even subcalculus. Consequently, 
\begin{equation}
\begin{split}
&\iota_D^*\tr (\A_{\mu}^k e^{-t\Delta_k(\mu)}) \sim \rho_{\ff}^{-1-b}
\left(\sum\limits_{j=0}^{\infty}a_{2j}\rho_{\ff}^{2j}\right), \\
&\iota_D^*\tr (\A_{\mu}^k e^{-t\Delta_k(\mu)}) \sim \rho_{\td}^{-m}
\left(\sum\limits_{j=0}^{\infty}b_{2j}\rho_{\td}^{2j}\right).
\end{split}\end{equation}
By the pushforward theorem, we now have
\begin{align}
\Tr \left(\A^k_{\mu}e^{-t\Delta_k(\mu)}\right) \sim \sum_{j=0}^{\infty}A_jt^{j-\frac{m}{2}}+\sum_{j=0}^{\infty}C_jt^{j-\frac{b}{2}}+
\sum_{j\in \mathfrak{I}}G_jt^{j-\frac{b}{2}}\log t,
\end{align}
where $\mathfrak{I}=\emptyset$ if $(m-b)$ is odd and $\mathfrak{I}=\N_0$ if $(m-b)$ is even.

Under either of these two sets of assumptions, we have shown that for each $k$, 
$$
\mathrm{Res}_{s=0}\int_0^{\infty}t^{s-1}\Tr \left(\A_{\mu}^k \left(e^{-t\Delta_k(\mu)}-
P_k(\mu)\right)\right)\, dt=-\Tr (\A_{\mu}^kP_k(\mu)).
$$
In view of the additional factor of $s$ in \eqref{s-factor}, we find
\begin{equation}
\frac{d}{d\mu}\log T(M,g_{\mu}) = \frac{1}{2}\sum_{k=0}^{m}(-1)^{k}\Tr \left(\A_{\mu}^k P_k(\mu)\right)
=\frac{d}{d\mu}\log \|\cdot \|^{-1}_{\det H^*(M,E), g_{\mu}}.
\end{equation}
This proves our central result.
\begin{thm}
Let $(M^m,g_{\mu})$ be a one-parameter family of asymptotically admissible edge metrics on a compact
space of odd dimension with simple edge $B^b$.  Assume that either
\begin{itemize}
\item[i)] each metric $g_\mu$ is strongly asymptotically admissible \\ with $g_\mu = g_{\mu_0} + h_\mu$ and
$|h_\mu|_{g_{\mu_0}} = \calO(x^{1+b})$, or
\item[ii)] each metric $g_\mu$ is even asymptotically admissible, and $b$ is odd.
\end{itemize}
Then the family of analytic torsion norms $\|\cdot \|^{RS}_{(M,g_{\mu})}$ is independent of the parameter $\mu$:
\[
\frac{d\, }{d\mu} \|\cdot \|^{RS}_{(M,g_{\mu})}=0.
\]
\end{thm}
We anticipate that this metric independence will have some interesting applications.

\section*{Appendix: The composition formul\ae}\label{composition}
In this Appendix we prove the main composition result for the incomplete edge heat calculus, and then show
that composition preserves the even subcalculus. This is in some sense the most technically difficult part of the parametrix
construction in the heat calculus. However, our method for understanding these compositions follows a now standard 
pattern in geometric microlocal analysis using Melrose's pushforward theorem for polyhomogeneous distributions 
with respect to $b$-fibrations between manifolds with corners.   We begin by reviewing this general result and then describe
how it implies the composition formula.  The composition formula in the heat calculus is almost identical to the one
in the conic setting proved by Mooers \cite{Moo}, and is also very closely related to one for the pseudodifferential
edge calculus in \cite{Maz-edge}; we refer also to the survey \cite{AM} for more on these matters. Because of all
these resources, we shall be brief here. 

To set the stage, however, suppose that $A \in \Psi^{k,\ell,\calE}_{\eh}(M)$ and $B \in \Psi^{k',\ell',\calE'}_{\eh}(M)$, with 
Schwartz kernels $K_A$ and $K_B$, which are pushforwards under the blowdown $\beta: M^2_h \to \RR^+ \times \wt{M}^2$ 
of polyhomogeneous distributions $\kappa_A$ and $\kappa_B$, respectively. We view them as sections of 
${}^{ie}\Lambda^p M \boxtimes {}^{ie}\Lambda^p M$ under rescaling $\Phi$. Let $C := A \circ B$, with Schwartz
kernel $K_C$ lifting to $\kappa_C$ on $M^2_h$. Our main task here is to show that $\kappa_C$ is polyhomogeneous,
with index sets at each face given as in the statement of the theorem in terms of those for $\kappa_A$ and $\kappa_B$,
respectively. 

To explain how this relates to the pushforward theorem, first observe that
\begin{equation}
K_{C}(t,x,y,z, x',y',z') = \int_0^t \int_{\wt{M}} K_A( t-t', w, w') K_B( t', w',  w'')\, dt' dV_g(w') 
\label{comp1}
\end{equation}
This can be rephrased as follows.  Consider the triple-space $\RR^+_{t'} \times \RR^+_{t''} \times \wt{M}_w \times \wt{M}_{w'} \times
\wt{M}_{w''}$, and the three projections
\begin{equation}
\begin{split}
&\pi_C :M^3\times \R^+_{t'}\times \R^+_{t''} \to M^2\times \R_{t'+t''}, \quad  (t',t'',w,w',w'') \to (w,w'',t'+t''), \\
&\pi_L: M^3\times \R^+_{t'}\times \R^+_{t''} \to M^2\times \R_{t'}, \qquad (t',t'',w,w',w'') \to (w,w',t'), \\
&\pi_R: M^3\times \R^+_{t'}\times \R^+_{t''} \to M^2\times \R_{t''}, \qquad (t',t'',w,w',w'') \to (w',w'',t'').
\end{split}
\label{projections}
\end{equation} 
Assuming that all kernels are extended to vanish for negative times, 
and reinterpreting them as densities in a suitable way specified below, we can rewrite \eqref{comp1} as
\[
K_C = (\pi_C)_* \left( \pi_L^* K_A  \pi_R^* K_B \right).
\]

The main idea is that we define a heat triple-space $M^3_h$; this will be obtained from $(\RR^+)^2 \times (\wt{M})^3$
by a sequence of blowups, and has the property that there are maps
\[
\Pi_L, \Pi_C, \Pi_R:  M^3_h \longrightarrow M^2_h
\]
which `cover' the three projections defined above. Lifting the composition to these spaces leads to the key formula
\begin{equation}
\kappa_C = (\Pi_C)_* \left( \Pi_L^* \kappa_A  \Pi_R^* \kappa_B \right).
\label{comp2}
\end{equation}
Thus it suffices to show that if $\kappa_A$ and $\kappa_B$ are polyhomogeneous, then their lifts to $M^3_h$
are also polyhomogeneous, and so too the product of these lifts, and most importantly, that the pushforward
by $\Pi_C$ of this product is again polyhomogeneous on $M^2_h$.  

\subsection*{The pushforward theorem}
We have already encountered a special case of this in the discussion preceding Theorem~\ref{heat-trace-corr}. First introduce 
some terminology.  Let $X$ and $X'$ be two compact manifolds with corners, and let $f: X \to X'$ 
be a smooth map.  Let $\{H_i\}$ and $\{H_j'\}$ be enumerations of the codimension one boundary faces of $X$ and $X'$,
respectively, and let $\rho_i$, $\rho_j'$ be global defining functions for $H_i$, resp.\ $H_j'$.  We say that the map $f$
is a $b$-map if $f^* \rho_j'$ is a smooth nonvanishing multiple of some product of nonnegative integer powers
of the defining functions $\rho_i$, or symbolically,
\[
f^* \rho_i' = A_{ij} \prod_j \rho_j^{e(i,j)}, \quad A_{ij} > 0,\ e(i,j) \in \mathbb{N} \cup \{0\}.
\]
This simply means that $f$ respects the boundary structure of these two spaces, and in particular maps each $H_i$ into some $H_j'$ 
with constant normal order of vanishing along the entire face. 

The map $f$ is called a $b$-submersion if $f_*$ induces a surjective map between the $b$-tangent bundles
of $X$ and $X'$. (The $b$-tangent space at a point $p$ of $\del X$ on a codimension  $k$ corner is spanned locally by  the
sections $x_1 \del_{x_1}, \ldots, x_k \del_{x_k}, \del_{y_j}$, where $x_1, \ldots, x_k$ are the defining functions for the faces
meeting at $p$ and the $y_j$ are local coordinates on the corner through $p$.)  If, in addition, the matrix $e(i,j)$ defined
above has the property that for each $j$ there is at most one $i$ such that $e(i,j) \neq 0$ (this condition simply means
that each hypersurface face $H_i$ in $X$ gets mapped into \emph{at most one} $H_j'$ in $X'$, or in other words, no
hypersurface in $X$ gets mapped to a corner in $X'$), then $f$ is called a $b$-fibration. 

We already introduced the notion of a $b$-density in \S 4.  Let us fix smooth nonvanishing $b$-densities $\nu_b$ on $X$
and $\nu_b'$ on $X'$. 
\begin{prop}[The Pushforward Theorem]
Let $f:X\to X'$ be a $b$-fibration. Let $u$ be a polyhomogeneous function on $X$ with index sets $E_j$ at the faces $H_j$ of $X$.  Suppose that each $(z,p) \in E_j$
has $\mbox{Re}\, z > 0$ if the index $j$ satisfies $e(i,j) = 0$ for all $j$ (which means that $H_j$ is mapped to the interior of $X'$).
Then the pushforward $f_* (u \nu_b)$ is well-defined and equals $h \nu_b'$ where $h$ is polyhomogeneous on $X'$ and has
an index family $f_b(\calE)$ given by an explicit formula in terms of the index family $\calE$ for $X$.
\end{prop}
For precise definition of the index family $f_b(\calE)$ see \cite{Mel2}.
Rather than giving the formula for the image index set in generality, let us describe it slightly informally but specifically
enough for the present situation.  If $H_{j_1}$ and $H_{j_2}$ are both mapped to a face $H_i'$, and if $H_{j_1} \cap H_{j_2} = \emptyset$,
then they contribute the index set $E_{j_1} \cup E_{j_2}$ to $H_i'$. If they do intersect, however, then the contribution is the
extended union $E_{j_1} \overline{\cup} E_{j_2}$. 

\subsection*{The triple space}
We now construct the reduced heat triple space $M^3_{\rh}$. The steps are dictated strictly by the requirements that the maps
$\pi_L, \pi_C, \pi_R$ all lift to $b$-fibrations $\Pi_L, \Pi_C, \Pi_R$, and that the space $M^3_{\rh}$ be as `small' as possible.
We revert to adapted boundary coordinates $w = (x,y,z)$.

The triple space $M^3_{\rh}$ is a parabolic blowup of $\wt{M}^3\times \R^+_{t'}\times \R^+_{t''}$, where the standard local 
coordinates in the singular edge neighborhoods of each of the three copies of $M$ are $(x,y,z)$, $(x',y',z')$ 
and $(x'',y'',z'')$. First blow up
\[
F=\{t'=t''=0, x=x'=x'' = 0, y=y'=y''\},
\]
parabolically with respect to both $t'$ and $t''$; then blow up the resulting space $[\wt{M}^3\times \R^+\times \R^+, F, \{dt', dt''\}]$ 
at the lift of $\mathcal{O}=\{t'=t''=0\}$; finally blow up the resulting space, parabolically in the respective time directions, at the 
lifts of each of the three submanifolds
\begin{equation}
\begin{split}
F_C&=\{t'=t''=0, x=x''=0, y=y''\}, \\
F_L&=\{t''=0, x'=x''=0, y'=y''\}, \\
F_R&=\{t'=0, x=x'=0, y=y'\}.
\end{split}
\end{equation} 
Thus altogether,
\[
M^3_{\rh} := \left[\wt{M}^3\times \R^+\times \R^+, F, \{dt'\}, \{dt''\}; \calO; F_C, \{dt', dt''\}; F_L, \{dt''\}; F_R, \{dt'\} \right].
\]

This is a rather difficult space to visualize, but one may `see' part of it by ignoring the time directions; 
the spatial part of $M^3_{\rh}$ is then exactly the same as the triple space appearing in the elliptic theory of edge operators, 
see \cite{Maz-edge}, which can be pictured as in Figure 5.
\begin{figure}[h]
\begin{center}
\begin{tikzpicture}

\draw  (-1,0) .. controls (-0.8,0) and (-0.5,-0.3) .. (-0.5,-0.5);
\draw  (1,0) .. controls (0.8,0) and (0.5,-0.3) .. (0.5,-0.5);
\draw  (-0.4,0.8) .. controls (0,0.6) and (0,0.6) .. (0.4,0.8);

\draw  (-1,0) -- (-2,-0.5);
\draw  (-0.5,-0.5) -- (-1.5,-1);

\draw  (1,0) -- (2,-0.5);
\draw  (0.5,-0.5) -- (1.5,-1);

\draw   (-0.4,0.8)  --  (-0.4,1.8) ;
\draw   (0.4,0.8)  --  (0.4,1.8);

\draw  (-1,0) .. controls (-1,0.2) and (-0.6,0.8) ..  (-0.4,0.8);
\draw  (1,0) .. controls (1,0.2) and (0.6,0.8) ..  (0.4,0.8);
\draw  (-0.5,-0.5) .. controls (-0.5,-0.6) and (0.5,-0.6) .. (0.5,-0.5);

\node at (0,0) {111};
\node at (-2,-0.9) {011};
\node at (2,-0.9) {110};
\node at (0,2) {101};
\node at (0,-1) {010};
\node at (-1.1,0.9) {001};
\node at (1.1,0.9) {100};
\end{tikzpicture}
\end{center}
\label{triple-space}
\caption{The spatial component of the triple space $M^3_{\rh}$.}
\end{figure}
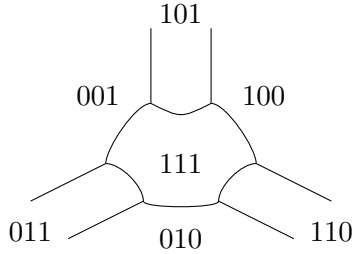 

Here, $(101), (011)$ and $(110)$ label the boundary faces created by blowing up $F_C, F_L$ and $F_R$, respectively. 
The face $(111)$ is the front face introduced by blowing up $F$. We denote the defining function for the face $(ijk)$  
by $\rho_{ijk}$. 

Now recall the projections $\pi_C$, $\pi_L$ and $\pi_R$ defined in \eqref{projections}.  These induce projections $\Pi_C$, $\Pi_L$
and $\Pi_R$ from 
$M^3_{\rh}$ to the reduced heat space $\wt{M}^2_{\rh}:=[M^2\times \R^+; FD; \{dt\}]$ introduced in \S 3 (see Figure 2). 
It is not hard to check that the choice of submanifolds that have been blown up ensures that these 
are in fact $b$-fibrations. 

Denote the defining functions for the right, front and left faces of each copy of $M^2_{\rh}$ by 
$\{\rho_{10}, \rho_{11}, \rho_{01}\}$, respectively.  These lift via the projections according to the following rules
\begin{equation}
\begin{split}
&\Pi_C^*(\rho_{ij})=\rho_{i0j}\rho_{i1j}, \\ 
&\Pi_L^*(\rho_{ij})=\rho_{ij0}\rho_{ij1}, \\ 
&\Pi_R^*(\rho_{ij})=\rho_{0ij}\rho_{1ij}.
\end{split}
\label{RLC}
\end{equation}

Now consider the behaviour in the time variables. Let $\beta^{(3)}: M^3_{\rh} \to (\RR^+)^2 \times \wt{M}^3$
be the blowdown map, $\tau_\calO$  the defining function for the face obtained by blowing up $\calO$,
and $\tau'$, $\tau''$ the defining functions for the two boundary faces in $M^3_{\rh}$ corresponding to $\{t'=0\}$ and $\{t''=0\}$. 
Then we find 
\begin{equation}
\label{beta-tau}
\begin{split}
&(\beta^{(3)})^*(t') = \tau'\tau_O\rho_{111}^2\rho_{110}^2\rho_{101}^2, \\
&(\beta^{(3)})^*(t'') = \tau''\tau_O\rho_{111}^2\rho_{011}^2\rho_{101}^2, \\
&(\beta^{(3)})^*(t'+t'') = \tau_O\rho_{111}^2\rho_{101}^2.
\end{split}
\end{equation} 
Let $\beta^{(2)}:M^2_{\rh}\to \R^+\times M^2$ be the blowdown map for the reduced heat space.
The lifts  $\left(\beta^{(2)}\right)^*(t'),\left(\beta^{(2)}\right)^*(t'')$ and 
$\left(\beta^{(2)}\right)^*(t'+t'')$ to $\Pi_L(M^3_{rh})$, $\Pi_R(M^3_{\rh})$ and $\Pi_C(M^3_{\rh})$, respectively, are 
equal to $T \rho_{11}^2$, where $T$ is the defining function for $\tf$ in $M^2_{\rh}$. Note 
\[
\Pi^*_{C,L,R} \circ \left(\beta^{(2)}\right)^* = \left(\beta^{(3)}\right)^* \circ \pi^*_{C,L,R}.
\]
Consequently, in view of \eqref{beta-tau} and \eqref{RLC}, we conclude
\begin{equation}
\begin{split}
&\Pi_C^*(T)=\tau_O, \\
&\Pi_L^*(T)=\tau'\tau_O\rho_{101}^2, \\
&\Pi_R^*(T)=\tau''\tau_O\rho_{101}^2.
\end{split}
\end{equation}

Using these data, we now derive the anticipated composition formula. Consider two elements of the heat calculus 
$A\in \Psi^{\ell,p,E_{\lf}, E_{\rf}}_{\eh}(M)$ and $B\in \Psi^{\ell',\infty, E'_{\lf}, E'_{\rf}}_{\eh}(M)$. We consider the Schwartz kernel of
each as a `right density', $K_A(t',x,y,z,x',y',z') dt' dx'dy'dz'$ and $K_B(t'',x',y',z',x'',y'',z'') dt'' dx'' dy'' dz''$. Then 
their product on $(\RR^+)^2 \times M^3$ is
\[
K_A(t',x,y,z,x',y',z') K_B(t'',x',y',z',x'',y'',z'') dt' dt'' dx' dy' dz' dx'' dy'' dz''.
\]
The integral over $dx' dy' dz'$ and over $t' + t'' = t$ gives 
$K_{A\circ B}(t,x,y,z,x'',y'',z'')dt \, dx''dy''dz''$. To put this into the same form required
in the pushforward theorem, multiply this expression by $dx dy dz$. 

Blowing up a submanifold of codimension $n$ amounts in local coordinates to introducing polar coordinates, so that 
the coordinate transformation of a density leads to $(n-1)^{\mathrm{st}}$ power of the radial function, which 
is the defining function of the corresponding front face. Hence we compute the lift
\begin{equation}
\begin{split}
&(\beta^{(3)})^* ( dt' dt'' dx dy dz dx' dy' dz' dx'' dy'' dz'') \\
&=\rho_{111}^{2b+6}\rho_{101}^{b+5}\rho_{110}^{b+3}\rho_{011}^{b+3}\tau_O \, 
\nu^{(3)}=\rho_{111}^{2b+6}\rho_{101}^{b+5}\rho_{110}^{b+3}\rho_{011}^{b+3}
\tau_O^2 \tau' \tau'' \left( \Pi \rho_{ijk} \right)
\nu^{(3)}_b,
\end{split}
\label{lift1}
\end{equation}
where $\nu^{(3)}$ is a density on $M^3_{\rh}$, smooth up to all boundary faces and everywhere 
nonvanishing; $\nu^{(3)}_b$ is a $b$-density, obtained from $\nu^{(3)}$ by dividing 
by a product of all defining functions on $M^3_{\rh}$; and $\left( \Pi \rho_{ijk} \right)$ is a product over all $(ijk)\in \{0,1\}^3$. 
Furthermore, the infinite order vanishing of $\kappa_A$ and $\kappa_B$ at $T=0$ implies that the product of the lifts
vanishes to infinite order in $\rho_{101}\tau_\calO \tau' \tau''$. Altogether, we obtain that
\[
\left(\Pi_L^*\kappa_A\right) \left(\Pi_R^*\kappa_B\right) (\beta^{(3)})^* ( dt' dt'' dx \, dy \, dz \, dx' dy' dz' dx'' dy'' dz'')  
 =\rho_{111}^{\ell + \ell'} \left( \Pi \rho_{ijk} \right) G \nu_b^{(3)},
\]
where $G$ is a polyhomogeneous function on $M^3_{\rh}$, vanishing to infinite order in $\rho_{101}\tau_\calO \tau' \tau''$, 
with index sets $E'_{\lf}$, $E_{\rf}$ and $E_{\lf} + E'_{\rf}$ 
at the faces $(001)$, $(100)$ and $(010)$, respectively. Moreover, $G$ has index sets $E_{\lf}+\ell'$ and 
$E'_{\rf}+\ell$ at the faces $(011)$ and $(110)$, respectively.

Note that since $\kappa_A$ does not vanish to infinite order on $\td$, the lift $\Pi_L^* \kappa_A$ is not polyhomogeneous
on $M^3_{\rh}$. Fortunately, the other factor $\kappa_B$ does vanish to infinite order there, and hence the product 
$\Pi_L^*\kappa_A \cdot \Pi_R^*\kappa_B$ is indeed polyhomogeneous on $M^3_{\rh}$.  Applying the Pushforward Theorem now gives
\begin{equation}
\begin{split}
&\left(\Pi_C\right)_*\left(\left(\Pi_L^*\kappa_A\right) \left(\Pi_R^*\kappa_B\right) (\beta^{(3)})^* ( dt' dt'' dx \, dy \, dz \, dx' dy' dz' dx'' dy'' dz'')\right) \\
& = \left(\beta^{(2)}\right)^*\left(K_{A\circ B}(t,x,y,z,x'',y'',z'')dt \, dx\, dy\, dz\, dx''dy''dz''\right) \\
&=  (\rho_{10}\rho_{01}\rho_{11} T) \, {\mathcal G} \, \nu^{(2)}_b,
\end{split}
\label{lift3}
\end{equation}
where $\nu^{(2)}_b$ is a $b$-density on $M^2_{\rh}$, and $\beta^{(2)}:M^2_{\rh} \to \R^+\times \wt{M}^2$ the corresponding 
blowdown map, and ${\mathcal G}$ is a polyhomogeneous function on $M^2_{\rh}$, which vanishes to infinite order in $T$, 
has leading order $\ell + \ell'$ without additional log-terms in its asymptotic behaviour at $\ff$ due to 
infinite order vanishing at $(101)$, and with the index sets $(P_{\lf}, P_{\rf})$ at the left and right boundary faces 
arising as extended unions (recall \eqref{extended})
\begin{equation}
\begin{split}
P_{\lf}&=E'_{\lf}\overline{\cup} (E_{\lf}+ \ell'), \\
P_{\rf}&=E_{\rf}\overline{\cup} (E'_{\rf}+ \ell),
\end{split}
\end{equation} 
By an argument similar to \eqref{lift1}, we compute
\begin{equation}
\begin{split}
\left(\beta^{(2)}\right)^*(dt \, dx\, dy\, dz\, dx''dy''dz'')=\rho_{11}^{b+3}\left(\rho_{10}\rho_{11}\rho_{01} T\right) \nu^{(2)}_b .
\end{split}
\label{lift4}
\end{equation}
Consequently, combining \eqref{lift3} and \eqref{lift4}, we deduce that $(\beta^{(2)})^* K_{A \circ B}=\kappa_{A\circ B}$ has index sets
$P_{\lf}$ and $P_{\rf}$ at the faces $(01)$ and $(10)$ of $M^2_{\rh}$, is $\rho^{-b-3 + \ell +\ell'}$ times a smooth
function at $\ff$, and vanishes to infinite order in $T$. 

This proves 
\begin{thm}
For index sets $E_{\lf}$ and $E'_{\rf}$ such that $E_{\lf}+E'_{\rf}>-1$, we have
$$\Psi^{l,p,E_{\lf}, E_{\rf}}_{e-h}(M) \circ 
\Psi^{l',\infty,E'_{\lf}, E'_{\rf}}_{e-h}(M) \subset 
\Psi^{l+l',\infty,P_{\lf}, P_{\rf}}_{e-h}(M),$$
where the front face expansion does not contain logarithmic terms and 
\begin{equation*}
\begin{split}
P_{\lf}&=E'_{\lf}\overline{\cup} (E_{\lf}+ \ell'), \\
P_{\rf}&=E_{\rf}\overline{\cup} (E'_{\rf}+ \ell). 
\end{split}
\end{equation*} 
\end{thm}

\section*{Composition in the even subcalculus}
We now give the proof of Proposition~\ref{evencomp}, and show that composition preserves the subspace of 
operators which satisfy the parity conditions from Definition~\ref{evencalcdef}.  In fact, we refer to \cite{Mel-APS}
for a proof that operators which satisfy these conditions near $\td$ are closed under composition, since the
proof is exactly the same here. Thus we focus on behaviour near $\ff$.  It clearly suffices to consider operators
with Schwartz kernels supported near some compact subset of the interior of the front face. Hence we are free
to write out and examine the composition using projective coordinates.

Note first that the pushforward by $\Pi_C$ does not introduce logarithmic 
terms in the front face expansion of $\kappa_{A\circ B}$, since 
the kernel on $M^3_{\rh}$ is vanishing to infinite order at $(101)$. 
Hence, for $\kappa_A$ and $\kappa_B$ with integer exponents in their front face expansions, 
$\kappa_{A\circ B}$ has a front face expansion with integer exponents as well, and 
we may check its evenness according to Definition~\ref{evencalcdef}.

As in the `full' composition formula, we lift the integrand in the formula 
\[ 
K_{A \circ B}(t,w) = \int_0^t \int_M K_A(t-t',w,w') K_B(t',w',\wt{w}) \, dt' dw'
\]
to $M^3_{\rh}$.  For the left factor, use the coordinates \eqref{top-coord}, while on the right we
introduce a further singular coordinate change
\[
v=\frac{y'-y''}{\rho}, \ \sigma=\frac{t'}{t}, \xi'=\frac{x'}{\rho}, 
\]
where $\rho = \sqrt{t}$. Following the substitutions in formula (18) of \cite{Gr}, we obtain
\begin{multline*}
\kappa_{A \circ B} (\rho, \xi, \widetilde{\xi}, w, \widetilde{y}, z, \widetilde{z})= \\ 
\rho^{b+3}\int_0^1\int \kappa_A (\rho\sqrt{1-\sigma},\frac{\xi}{\sqrt{1-\sigma}}, 
\frac{\xi'}{\sqrt{1-\sigma}}, \frac{w-v}{\sqrt{1-\sigma}}, \widetilde{y}+v\rho) \\ 
\kappa_B (\rho\sqrt{\sigma},\frac{\xi'}{\sqrt{\sigma}}, \frac{\widetilde{\xi}}{\sqrt{\sigma}}, 
\frac{v}{\sqrt{\sigma}}, \widetilde{y}) \, d\sigma \, d\xi' \, dv \, dz'. 
\end{multline*}

Now expand the two factors to get:
\begin{multline*}
\kappa_A (\rho\sqrt{1-\sigma},\frac{\xi}{\sqrt{1-\sigma}}, \frac{\xi'}{\sqrt{1-\sigma}}, \frac{w-v}{\sqrt{1-\sigma}}, 
\widetilde{y}+v\rho) \\ 
\sim (\rho\sqrt{1-\sigma})^{-b-3}\sum_k A_k(\frac{\xi}{\sqrt{1-\sigma}}, \frac{\xi'}{\sqrt{1-\sigma}}, 
\frac{w-v}{\sqrt{1-\sigma}}, \widetilde{y}+v\rho)(\rho\sqrt{1-\sigma})^k,
\end{multline*}
\[
\kappa_B(\rho\sqrt{\sigma},\frac{\xi'}{\sqrt{\sigma}}, \frac{\widetilde{\xi}}{\sqrt{\sigma}}, \frac{v}{\sqrt{\sigma}}, 
\widetilde{y}) \sim (\rho\sqrt{\sigma})^{-b-3}
\sum_k B_k(\frac{\xi'}{\sqrt{\sigma}}, \frac{\widetilde{\xi}}{\sqrt{\sigma}}, \frac{v}{\sqrt{\sigma}}, \widetilde{y}) (\rho\sqrt{\sigma})^k.
\]
By assumption, $A_{2k}$ and $B_{2k}$ are even while $A_{2k+1}$ and $B_{2k+1}$ are odd in their third slots.

Since each $A_k$ is smooth in the final variable, we expand further to obtain
\[
A_k(\frac{\xi}{\sqrt{1-\sigma}}, \frac{\xi'}{\sqrt{1-\sigma}}, \frac{w-v}{\sqrt{1-\sigma}}, \widetilde{y}+v\rho)\sim
\sum_{j=0}^{\infty}A_{(k,j)}(\frac{\xi}{\sqrt{1-\sigma}}, \frac{\xi'}{\sqrt{1-\sigma}}, \frac{w-v}{\sqrt{1-\sigma}}, \widetilde{y},v)\rho^j.
\]
Note that clearly, 
\[
A_{(k,j)}(\frac{\xi}{\sqrt{1-\sigma}}, \frac{\xi'}{\sqrt{1-\sigma}}, \frac{w-v}{\sqrt{1-\sigma}}, \widetilde{y},-v)=
(-1)^jA_{(k,j)}(\frac{\xi}{\sqrt{1-\sigma}}, \frac{\xi'}{\sqrt{1-\sigma}}, \frac{w-v}{\sqrt{1-\sigma}}, \widetilde{y},v).
\]

Now define
\[
A_{[i]}:=\sum\limits_{k+j=i}A_{(k,j)}.
\]
Then the two parity conditions together give that
\[
A_{[i]}(\frac{\xi}{\sqrt{1-\sigma}}, \frac{\xi'}{\sqrt{1-\sigma}}, -\frac{w-v}{\sqrt{1-\sigma}}, 
\widetilde{y},-v)=(-1)^i A_{[i]}(\frac{\xi}{\sqrt{1-\sigma}}, \frac{\xi'}{\sqrt{1-\sigma}}, \frac{w-v}{\sqrt{1-\sigma}}, 
\widetilde{y},v).
\]

Inserting these expansions into the integral leads to an expansion of $\kappa_{A \circ B}$.  Since the rescaled volume form $x^{-f}dV_g$ is  a
smooth function of $\rho v$ and $\rho^2$, it does not affect these parity considerations. 

Collecting all the factors involving only $\sigma$ into a single function $f(\sigma)$, we get
\begin{multline*}
\kappa_{A \circ B}(\rho, \xi, \widetilde{\xi}, w, \widetilde{y}, z, \widetilde{z}) \sim  
\rho^{-b-3}\sum_{i,k} \rho^{i+k}\int_0^1\int f(\sigma) A_{[i]}(\frac{\xi}{\sqrt{1-\sigma}}, \frac{\xi'}{\sqrt{1-\sigma}}, 
\frac{w-v}{\sqrt{1-\sigma}}, \widetilde{y},v) \\
\times B_k(\frac{\xi'}{\sqrt{\sigma}}, \frac{\widetilde{\xi}}{\sqrt{\sigma}}, \frac{v}{\sqrt{\sigma}}, \widetilde{y}) \,
d\sigma \, d\xi' \, dv\, dz' =:  \rho^{-b-3}\sum_{i,k} \rho^{i+k}(A\circ B)_{i,k}(\xi, \widetilde{\xi}, w, \widetilde{y}, z, \widetilde{z}).
\end{multline*}

Finally, observe that 
\begin{align*}
&(A\circ B)_{i,k}(\xi, \widetilde{\xi}, -w, \widetilde{y}, z, \widetilde{z})\\ 
&=\int_0^1\int f(\sigma) A_{[i]}(\frac{\xi}{\sqrt{1-\sigma}}, \frac{\xi'}{\sqrt{1-\sigma}}, 
-\frac{w+v}{\sqrt{1-\sigma}}, \widetilde{y},v) B_k(\frac{\xi'}{\sqrt{\sigma}}, \frac{\widetilde{\xi}}{\sqrt{\sigma}}, 
\frac{v}{\sqrt{\sigma}}, \widetilde{y})\\
&= (-1)^i\int_0^1\int f(\sigma) A_{[i]}(\frac{\xi}{\sqrt{1-\sigma}}, \frac{\xi'}{\sqrt{1-\sigma}}, 
\frac{w+v}{\sqrt{1-\sigma}}, \widetilde{y},-v) B_k(\frac{\xi'}{\sqrt{\sigma}}, \frac{\widetilde{\xi}}{\sqrt{\sigma}}, 
\frac{v}{\sqrt{\sigma}}, \widetilde{y})\\
&=(-1)^{i+k}\int_0^1\int f(\sigma) A_{[i]}(\frac{\xi}{\sqrt{1-\sigma}}, \frac{\xi'}{\sqrt{1-\sigma}}, 
\frac{w+v}{\sqrt{1-\sigma}}, \widetilde{y},-v) B_k(\frac{\xi'}{\sqrt{\sigma}}, \frac{\widetilde{\xi}}{\sqrt{\sigma}}, -\frac{v}{\sqrt{\sigma}}, \widetilde{y})\\
&=(-1)^{i+k}\int_0^1\int f(\sigma) A_{[i]}(\frac{\xi}{\sqrt{1-\sigma}}, \frac{\xi'}{\sqrt{1-\sigma}},
\frac{w-v}{\sqrt{1-\sigma}}, \widetilde{y},v) B_k(\frac{\xi'}{\sqrt{\sigma}}, \frac{\widetilde{\xi}}{\sqrt{\sigma}}, 
\frac{v}{\sqrt{\sigma}}, \widetilde{y})\\
&=(-1)^{i+k}(A\circ B)_{i,k}(\xi, \widetilde{\xi}, w, \widetilde{y}, z, \widetilde{z}).
\end{align*}
Hence the Schwartz kernel of $A\circ B$ satisfies the parity condition too. 

\bibliographystyle{amsplain}

\begin{thebibliography}{10}

\bibitem{Albin} P. Albin, 
\emph{A renormalized index theorem for some complete asymptotically regular metrics: The Gauss-Bonnet theorem}, 
Adv. in Math. {\bf 213} (2007), No. 1, 1-52.
\bibitem{ALMP} P. Albin, E. Leichtnam, R. Mazzeo, P. Piazza \emph{The signature package on Witt spaces, I. Index classes}, 
to appear, Ann. de la Ec. Norm. Sup. 
\bibitem{AM} P. Albin, R. Mazzeo \emph{Geometric constructions of heat kernels: a user's guide}, to appear.
\bibitem{AR} P. Albin, F. Rochon
\emph{Families index for manifolds with hyperbolic cusp singularities}, Int. Math. Res. Not. (2009), no. 4, 625-697.
\bibitem{BDV} E. Bahuaud, E.Dryden and B. Vertman \emph{Mapping properties of the Heat operator on edge manifolds}, 
arXiv:1105.5119. 
\bibitem{BGV} M. Berline, E. Getzler and M. Vergne \emph{Heat kernels and Dirac operators}, 
Grundlehren der Mathematischen Wissenschaften {\bf 298}, Springer-Verlag, Berlin (1992)
\bibitem{BS} J. Br\"uning, R. Seeley \emph{The resolvent expansion of second order regular singular operators}, J. Funct. Anal. {\bf 73} (1987), 369-429.
\bibitem{BS3} J. Br\"uning, R. Seeley \emph{The expansion of the resolvent near a singular stratum of conical type} J. Funct. Anal. {\bf 95} (1991), 255-290 
\bibitem{Ch1} J. Cheeger \emph{On the spectral geometry of spaces with conical singularities}, Proc. Nat. Acad. Sci. USA {\bf 74} (1979), 2651-2654 
\bibitem{Ch2} J. Cheeger \emph{Spectral geometry of singular Riemannian spaces}, J. Diff. Geom. {\bf 18} (1983), 575-657. 
\bibitem{Ch-AT} J. Cheeger \emph{Analytic Torsion and the Heat Equation},  Ann. of Math.(2) {\bf109} (1979) no. 2, 259--322.
\bibitem{DH} X. Dai, X. Huang, \emph{The intersection R-torsion of a finite cone}, Preprint (2010).
\bibitem{Dar} A. Dar \emph{Intersection R-torsion and analytic torsion for pseudo-manifolds}, Math. Z. {\bf 194} (1987), 193-216.
\bibitem{EMM} C.L. Epstein, R. Melrose and G. Mendoza 
\emph{Resolvent of the Laplacian on strictly pseudoconvex domains}, Acta Math. {\bf 167} (1991), no. 1-2, 1--106. 
\bibitem{GKM}  J. Gil, T. Krainer and G. Mendoza \emph{On the closure of elliptic wedge operators}, arXiv:1007.2397v2 [math.AP]
\bibitem{GM} J. Gil and G. Mendoza \emph{Adjoints of elliptic cone operators} Amer. J. Math. {\bf 125} (2003), no. 2, 357–408.
\bibitem{Gui} C. Guillarmou \emph{Meromorphic properties of the resolvent on asymptotically hyperbolic manifolds}, to appear,
Duke Math. Jour. 
\bibitem{Gr} D. Grieser \emph{Notes on heat kernel asymptotics}. Preprint, available at 
http://www.staff.uni-oldenburg.de/daniel.grieser/wwwvortraege/vortraege.html
\bibitem{HM} E. Hunsicker, R. Mazzeo 
\emph{Harmonic forms on manifolds with edges}, Int. Math. Res. Not. (2005), no. {\bf 52}, 3229--3272
\bibitem{Ha} A. Hassell \emph{Analytic surgery and analytic torsion} Comm. Anal. Geom. {\bf 6} (1998), no. 2, 255-–289.
\bibitem{KLP} K. Kirsten, P. Loya, J. Park, with an Appendix by Boris Vertman \emph{Exotic expansions and pathological properties of 
zeta-functions on conic manifolds}, J. Geom. Anal. {\bf 18} (2008), 835-888.
\bibitem{Lesch} M. Lesch \emph{Operators of Fuchs type, conical singularities and asymptotic methods}, Teubner Texte zur Mathematik Vol. 136, 
Teubner--Verlag, Leipzig, 1997.  Also available as arXiv:dg-ga/9607005. 
\bibitem{Maz-edge} 
R. Mazzeo \emph{Elliptic theory of differential edge operators, I}, Comm. Part. Diff. Eq. {\bf 16} (1991), No. 10, 1615-1664.
\bibitem{Mel-APS} 
R. Melrose \emph{The Atiyah-Patodi-Singer index theorem} Research Notes in Math., Vol. 4, A K Peters, Massachusetts (1993)
\bibitem{Mel2} 
R. Melrose \emph{Calculus of conormal distributions on manifolds with corners} Intl. Math. Research Notices, No. 3  (1992), 51-61.
\bibitem{Moo} 
E. Mooers \emph{Heat kernel asymptotics on manifolds with conic singularities}, J. Anal. Math. {\bf 78} (1999), 1-36.
\bibitem{Mue-AT} 
W. M\"uller \emph{Analytic Torsion and R-torsion of Riemannian manifolds}, Adv. Math. {\bf 28} (1978), 233-305.
\bibitem{RS} 
D.B. Ray, I.M. Singer \emph{R-torsion and the Laplacian on Riemannian manifolds}, Adv. Math. {\bf 7} (1971), 145-210.
\bibitem{Sch1} 
B. -W. Schulze \emph{Pseudo-differential operators on manifolds with singularities}, North-Holland, Amsterdam (1991)
\bibitem{Scott}
S. Scott
\emph{Traces and determinants of pseudodifferential operators} Oxford Mathematical Monographs. Oxford University Press, Oxford (2010).
\bibitem{Seeley}
R. Seeley \emph{Trace expansions for the Zaremba problem}, Comm. Part. Diff. Eq. {\bf 28} (2003), no. 3-4, 601--616.
\bibitem{Spr}
M. Spreafico \emph{The analytic torsion of a cone over a sphere}, J. Math. Pures Appl. (9) {\bf 93} (2010), no. 4, 408--435. 
\bibitem{Vert} 
B. Vertman \emph{Analytic torsion of a bounded generalized cone}, Comm. Math. Phys. {\bf 290} (2009), no. 3, 813--860. 
\bibitem{Zh}
W. Zhang  \emph{Lectures on Chern-Weil theory and Witten deformations} 
Nankai Tracts in Mathematics, 4. World Scientific Publishing Co., Inc., River Edge, NJ, 2001

\end{thebibliography}
\def\cprime{$'$} \def\polhk#1{\setbox0=\hbox{#1}{\ooalign{\hidewidth
  \lower1.5ex\hbox{`}\hidewidth\crcr\unhbox0}}}
\providecommand{\bysame}{\leavevmode\hbox
to3em{\hrulefill}\thinspace}
\providecommand{\MR}{\relax\ifhmode\unskip\space\fi MR }
\providecommand{\MRhref}[2]{%
  \href{http://www.ams.org/mathscinet-getitem?mr=#1}{#2}
} \providecommand{\href}[2]{#2}

\end{document}